\documentclass[class=article, final, crop=false, float = true]{standalone}

\usepackage{authblk}
\usepackage{setspace}
\usepackage{tiaArticleStyle}
\usepackage{bbold}
\usepackage{makecell}

\graphicspath{{media/}}
\usepackage{wrapfig, framed, caption}
\bibliography{references.bib}

\begin{document}

\title{Learning Regularization Parameters of Inverse Problems \\ via Deep Neural Networks}

\author[$\dagger$]{Babak Maboudi Afkham}
\author[$\star$]{Julianne Chung}
\author[$\star$]{Matthias Chung}

\affil[$\dagger$]{DTU Compute, Department of Applied Mathematics and Computer Science, Technical University of Denmark, Lyngby, Denmark}
\affil[$\star$]{Department of Mathematics, Virginia Tech, Blacksburg, VA 24060 USA}

\maketitle

\begin{abstract}
In this work, we describe a new approach that uses deep neural networks (DNN) to obtain regularization parameters for solving inverse problems.  We consider a supervised learning approach, where a network is trained to approximate the mapping from observation data to regularization parameters.  Once the network is trained, regularization parameters for newly obtained data can be computed by efficient forward propagation of the DNN. We show that a wide variety of regularization functionals, forward models, and noise models may be considered. The network-obtained regularization parameters can be computed more efficiently and may even lead to more accurate solutions compared to existing regularization parameter selection methods.  We emphasize that the key advantage of using DNNs for learning regularization parameters, compared to previous works on learning via optimal experimental design or empirical Bayes risk minimization, is greater generalizability.  That is, rather than computing one set of parameters that is optimal with respect to one particular design objective, DNN-computed regularization parameters are tailored to the specific features or properties of the newly observed data. Thus, our approach may better handle cases where the observation is not a close representation of the training set.  Furthermore, we avoid the need for expensive and challenging bilevel optimization methods as utilized in other existing training approaches.  Numerical results demonstrate the potential of using DNNs to learn regularization parameters.
\end{abstract}


Date: \today

\section{Introduction \& Background}
Many scientific problems can be modeled as
\begin{equation} \label{eq:invprob}
  \bfb = \bfA(\bfx_{\rm true}) + \bfvarepsilon,
\end{equation}
where $\bfx_{\rm true} \in \bbR^n$ is a desired solution, $\bfA:\bbR^n \to\bbR^m$ models some forward process mapping onto observations $\bfb \in \bbR^m$ at pre-determined design points, with unknown additive noise $\bfvarepsilon \in \bbR^m$. The goal in \emph{inverse problems} is to obtain an approximate solution $\widehat\bfx$ to $\bfx_{\rm true}$, given $\bfb$ and $\bfA(\cdot)$. However, solving inverse problems may be challenging due to ill-posedness, whereby a solution does not exist, is not unique, or does not depend continuously on the data \cite{engl1996regularization,hansen2010discrete}.  Regularization in the form of prior knowledge on the distribution of $\bfx_{\rm true}$ must be included to compute reasonable solutions.  There are many forms of regularization, and we consider variational regularization and regularization via early stopping techniques.  The goal of variational regularization is to minimize some \emph{loss function},
\begin{equation}\label{eq:loss}
  \min_{\bfx} \ \calJ(\bfx, \bfb) +  \calR(\bfx, \bflambda),
\end{equation}
where $\calJ:\bbR^n \times \bbR^m \to \bbR$ characterizes a \emph{model-data misfit} measuring the discrepancy between a model prediction and the observations $\bfb$ and the functional $\calR:\bbR^n \times \bbR^\ell \to \bbR$ represents a \emph{regularization} term defined by some parameters $\bflambda$. A commonly used model-data misfit is the (squared) Euclidean distance, i.e., $\calJ(\bfx, \bfb) = \norm[2]{\bfA(\bfx) - \bfb}^2$. We assume that the regularization term carries prior knowledge of the desired solution $\bfx_{\rm true}$ and that the parameters in $\bflambda$ define the regularity of the desired parameters in $\bfx$ and hence the regularization term. For instance, $\bflambda$ may contain one regularization parameter that determines the weight or strength of the regularization term, e.g., $\calR(\bfx,\lambda) = \lambda^2 \norm[2]{\bfx}^2$ corresponds to standard Tikhonov regularization. Another example arises in the identification of inclusions (e.g., cancers or other anomalies) in images, where $\lambda$ characterizes the regularity of the inclusion and must be estimated.  In other scenarios, $\bflambda$ may contain a set of parameters (e.g., for the prior covariance kernel function) that fully determine the regularization functional $\calR(\mdot,\bflambda)$.
For simplicity and illustrative purposes, we assume that optimization problem~\eqref{eq:loss} is sufficiently smooth, convex, and has a unique global minimizer $\widehat\bfx(\bflambda)$ for any suitable $\bflambda$.

A major computational difficulty in the solution of~\eqref{eq:loss} is that $\bflambda$ must be determined \textit{prior} to solution computation.  Selecting appropriate parameters $\bflambda$ can be a very delicate and computationally expensive task, especially for large-scale and nonlinear problems \cite{haber2000gcv,mead2008newton,vogel1996non,galatsanos1992methods,farquharson2004comparison}.  Common approaches for estimating the regularization parameters require solving \eqref{eq:loss} multiple times for various parameter choices, which may require solving many large-scale nonlinear optimization problems, until some criterion is satisfied.  For example, the discrepancy principle seeks parameters $\bflambda$ such that $\calJ(\widehat\bfx(\bflambda), \bfb) \approx T$ where $T$ is some target misfit (e.g., based on the noise level of the problem). Since the parameters $\bflambda$ determine the prior, they may also be referred to as hyper-parameters, and hierarchical prior models may be incorporated in a Bayesian formulation to include probabilistic information about the hyper-priors \cite{bardsley2018computational,tenorio2017introduction}.

For applications where training data is readily available or can be experimentally generated (e.g., via Monte Carlo simulations), supervised learning approaches have been used to learn regularization parameters or more generally ``optimal'' regularizers for inverse problems.
A new paradigm of obtaining regularizers was first introduced in~\cite{haber2003learning}.  In their groundbreaking, but often overlooked publication, Haber and Tenorio proposed a supervised learning approach to learn optimal regularizers. This framework leads to a bilevel optimization problem, where the inner problem consists of the underlying inverse problem assuming a fixed regularization functional. The outer problem -- often referred to as the design problem -- seeks an optimal regularization functional, given training data.
More specifically, given training data of true solutions and corresponding observations, $\left\{\bfx_{\rm true}^{(j)}, \bfb^{(j)}\right\}_{j=1}^J$ optimal regularization parameters are computed as
\begin{equation}
\label{eq:OED}
\min_\bflambda \ \tfrac{1}{2J}\sum_{j=1}^J \norm[2]{\widehat \bfx^{(j)}(\bflambda) - \bfx_{\rm true}^{(j)}}^2
\qquad
\mbox{ with } \qquad
\widehat\bfx^{(j)} (\bflambda) = \argmin_\bfx \ \calJ(\bfx, \bfb^{(j)}) + \calR(\bfx,\bflambda).
\end{equation}
This learning approach has shown great success in a variety of problems and has given rise to various new approaches for optimal experimental design \cite{haber2009numerical,haber2008numerical, chung2011designing,calatroni2017bilevel,ruthotto2018optimal} and for obtaining optimal regularizers \cite{chung2011designing,chung2012optimal,de2017bilevel}. Various other research groups build on the same bilevel supervised learning principle, e.g., \cite{de2016machine, antil2020bilevel} and references therein. A main challenge of this approach is to numerically solve the bilevel optimization problem.
We emphasize that computed parameters are expected to be optimal \emph{on average} or with respect to other design criteria and may fail in practice if the observation is very different than the training set, see~\cite{Atkinson2007}.

There is another class of supervised learning methods that has gained increased attention for solving inverse problems in recent years.  These methods exploit deep neural network (DNN) learning techniques such as convolution neural networks or residual neural networks \cite{lucas2018using,mccann2017convolutional}. Initially, these machine learning techniques were used for post-processing solutions, e.g., to improve solution quality or to perform tasks such as image classification \cite{zhang2017learning}. However, deep learning techniques have also been used for solving inverse problems.  The prevalent approach, especially in image processing, is to take an end-to-end approach or to use deep learning methods to replace a specific task (e.g., image denoising or deblurring). For example, in \cite{li2020nett} neural networks are used to learn the entire mapping from the data space to the inverse solution and in \cite{de2016machine,wang2020learning,hammernik2018learning} DNNs were used to learn the entire regularization functional.
Note, these approaches do not include domain-specific knowledge, but rather replace the inversion of a physical system with a black-box forward propagating process also referred to as surrogate modeling. Hence, the limitations of these approaches appear in the sensitivity of the network (e.g., to large dimensional input-output maps as they appear in imaging applications). Work on unsupervised learning approaches such as deep image priors have been considered as an alternative \cite{dittmer2020regularization}.
Another remedy is to reduce the size of the network inputs. In \cite{liu2021machine} a machine learning-based prediction framework is used to estimate the regularization parameter for seismic inverse problems.  Using a list of $19$ predefined features (e.g, including energy power and distribution characteristics of the data and residual) for the synthetic observation data, the authors use a random forest algorithm to train a decision tree for the task of regression.  Although the idea to learn the regularization parameter (representing the strength of regularization) is a special case of the framework we consider, the main distinction of our approach is that we consider DNNs to represent the mapping from the observation to the optimal regularization parameter.

In this work, we described a new approach to learn the parameters $\bflambda$ that define the regularization by training a neural network to learn a mapping from observation to regularization parameters. We begin by assuming that there exists a nonlinear \emph{target function} $\bfPhi: \bbR^m \to \bbR^p$ that maps an input vector $\bfb \in \bbR^m$ to a vector $\bflambda \in \bbR^p$,
\begin{equation}
\label{eq:phistar}
\bflambda = \bfPhi(\bfb).
\end{equation}
The function $\bfPhi$ is a nonlinear mapping that takes any vector in $\bbR^m$ (e.g., the observations) to a set of parameters in $\bflambda$ (e.g., the regularization parameters).  Thus, in the inverse problems context, we refer to $\bfPhi$ as an \emph{observation-to-regularization} mapping, and we assume that this function is well-defined.

A major goal of this work is to estimate the \emph{observation-to-regularization} mapping $\bfPhi$ by approximating it with a neural network and learning the parameters of the network.  We consider DNNs, which have gained increased popularity and utility in recent years due to their universal approximation properties \cite{cybenko1989approximation}.
That is, we assume that the observation-to-regularization mapping can be approximated using a feedforward network that is defined by some parameters $\bftheta$.
The network is a mapping $\widehat \bfPhi(\ \cdot\ ;\ \bftheta):\bbR^m \to \bbR^p$ that is defined by the weights and biases contained in $\bftheta$.  Given an \emph{input} $\bfb$, the \emph{output} of the network is given by
\begin{equation}
\label{eq:phi}
\widehat \bflambda(\bftheta) = \widehat \bfPhi(\bfb;\ \bftheta),
\end{equation}
see Figure~\ref{fig:nn} for a general schematic and Section~\ref{sec:NNlearning} for details of the network.
Notice that for a well-chosen set of parameters $\bftheta,$ the DNN can approximate the desired mapping, $\widehat \bfPhi \approx \bfPhi$, but a robust learning approach is needed to estimate network parameters $\bftheta$ that result in a good network approximation of the function.  More specifically, the goal is to minimize the Bayes risk, i.e., the expected value of some loss function $D:\bbR^p \times \bbR^p \to \bbR$. Let $\bfb = \bfA(\bfx_{\rm true}) + \bfvarepsilon$ where $\bfx_{\rm true}$ and $\bfvarepsilon$ are random variables. The learning problem can be written as an optimization problem of the form,
\begin{equation}
\label{eq:stochopt}
    \arg\min_\bftheta \ \bbE_{\bfx_{\rm true}, \bfvarepsilon} \ D \left(\widehat \bfPhi (\bfA(\bfx_{\rm true})+\bfvarepsilon;\ \bftheta ), \bflambda_{\rm opt} \right),
\end{equation}
where $\bflambda_{\rm opt}$ may be provided or computed. For example, for some problems, $\bflambda_{\rm opt}$ could be obtained by solving bilevel optimization problem,
\begin{equation}
\label{eq:lambdatrue}
\bflambda_{\rm opt} = \argmin_\bflambda \ \norm[2]{\widehat \bfx(\bflambda) - \bfx_{\rm true}}
\quad
\mbox{ with } \quad
\widehat\bfx (\bflambda) = \argmin_\bfx \ \calJ(\bfx, \bfb) + \calR(\bfx,\bflambda).
\end{equation}
The Bayes risk minimization problem \eqref{eq:stochopt} is a stochastic optimization problem, and the literature on stochastic programming methodologies is vast \cite{shapiro2014lectures}.  The learning problem can also be interpreted as an optimal experimental design (OED) problem, where the goal is to design a network to represent the observation-to-regularization mapping. This is different than OED problems that seek to optimize for the regularization parameters directly, e.g., \cite{haber2003learning,chung2011designing,haber2008numerical}.

Notice that the expected value in \eqref{eq:stochopt} is defined in terms of the distributions of $\bfx_{\rm true}$ and $\bfvarepsilon$, and thus if such knowledge is available or can be well-approximated, then Bayes risk minimization can be used.  However, for problems where the distributions of $\bfx_{\rm true}$ and $\bfvarepsilon$ are unknown or not obtainable, we consider empirical Bayes risk design problems, where training data or calibration data are used to approximate the expected value.  Assume that we have training data $\left\{\bfx_{\rm true}^{(j)}, \bfvarepsilon^{(j)}\right\}_{j = 1}^J$ and that the goal is to estimate the regularization parameters $\bflambda$ that are deemed optimal.  Then for each training sample, we would first obtain $\bflambda_{\rm opt}^{(j)}$ for $j = 1, \dots, J$ by \eqref{eq:lambdatrue}.  Then, we can approximate the Bayes risk problem \eqref{eq:stochopt} with the following empirical Bayes risk minimization problem,
\begin{equation}
    \label{eq:saa}
    \widehat \bftheta = \argmin_\bftheta \ \tfrac{1}{2J} \sum_{j=1}^J D\left(\widehat \bfPhi(\bfb^{(j)};\ \bftheta), \bflambda_{\rm opt}^{(j)}\right),
\end{equation}
where $\bfb^{(j)} = \bfA(\bfx_{\rm true}^{(j)})+ \bfvarepsilon^{(j)}.$
Given some loss function $D$, DNN $\widehat \bfPhi$, and the data set $\left\{\bfb^{(j)}, \bflambda_{\rm opt}^{(j)}\right\}_{j = 1}^J$, the goal of the supervised learning approach is to compute $\widehat \bftheta$ in an \emph{offline stage}.  Then in an \emph{online phase}, given a newly-obtained observation $\bfb$, regularization parameters for the new data can be easily and cheaply obtained via forward propagation through the network, i.e., $\widehat \bflambda = \widehat \bfPhi(\bfb; \ \widehat \bftheta)$. Once $\widehat \bflambda$ is computed and fixed, a wide range of efficient solution techniques may be used to solve the resulting inverse problem \eqref{eq:loss}.

\textbf{Overview of main contributions:}  In this work, we describe a new approach to estimate regularization parameters by training a DNN to approximate the mapping from observation data to regularization parameters.  Once the network is trained, regularization parameters may be computed efficiently via forward propagation through the DNN. There are various advantages of our proposed approach.
\begin{enumerate}
\item The DNN computed regularization parameters are tailored to the specific features or properties of the newly obtained data (e.g., the computed parameters are adapted to the amount of noise in the data). This potential for greater \emph{generalizability} to data that does not closely resemble features of the training set is a significant benefit compared to optimal experimental design or empirical Bayes risk minimization approaches where one set of design parameters is obtained that is (e.g., on average) good for the training set.
\item Given a new observation, the network-computed regularization parameters can be computed very \emph{efficiently} in an online phase, only requiring a forward propagation of the neural network.
Since this process requires only basic linear algebra operations and activation function evaluations, computing regularization parameters in this way is significantly faster than many existing regularization parameter selecting methods that may require solving multiple inverse problems for multiple parameter choices or may need derivative evaluations.
\item The DNN computed parameters can lead to solutions that are \emph{more accurate} than existing methods. Notice that after the regularization parameters are computed via a DNN, regularization is applied to the original problem and well-established solution techniques and software can be used to solve the resulting regularized problem. Contrary to black-box inversion methods, the physical forward model (which might be slightly different than the one used for training) is used during inversion.
\item A key advantage of the proposed work is the \emph{flexibility} of the approach in that a wide range of forward models and regularizers, most notably nonlinear ones, may be included, and the framework can learn other important features from data such as the degree of regularity of solutions.
\end{enumerate}
We also mention a few shortcomings of learning regularization parameters via DNNs. Certainly a downside of our method compared to full network inversion approaches is that we still require solving the resulting inverse problem in the online phase.  Another potential disadvantage is that as the dimension $p$ of the network output in \eqref{eq:phi} gets larger, more training data are required.

An overview of the paper is as follows.
Section \ref{sec:NNlearning} is dedicated to our proposed approach for learning a neural network for regularization parameter selection. We describe various details about the process from defining the network to optimization methods for learning. In Section~\ref{sec:oneparam} we focus on the special (and most common) case where we seek one regularization parameter corresponding to the strength of the regularization. Numerical experiments provided in Section \ref{sec:numerics} illustrate the benefits and potential of our approach for applications in tomography reconstruction, image deblurring, and diffusion.  Conclusions and future work are provided in Section \ref{sec:conclusions}.

\section{Parameter learning via training of neural networks}
\label{sec:NNlearning}
This section is devoted to our proposed approach to learn regularization parameters via DNNs for solving inverse problems.
We describe various components of our approach, but we begin with a general overview of the approach in Algorithm~\ref{alg:nn}.  Notice that there is an \emph{offline phase} and an \emph{online phase}.  In the offline phase, the training data is used to learn the network parameters. This requires solving a large scale optimization problem.  However, once the network parameters are computed, forward propagation of any new observation $\bfb$ through the network will produce a set of regularization parameters (e.g., for use in defining and solving the regularized problem).

\begin{algorithm}[bthp]
\setstretch{1.35}
    \begin{algorithmic}[1]
        \State \emph{offline phase}
        \State \hspace*{2ex} require model $\bfA(\mdot)$, noise model $\bfvarepsilon$, and $\bfx_{\rm true}^{(j)}$
        \State \hspace*{2ex} generate appropriate training signals $\bfb^{(j)} = \bfA(\bfx_{\rm true}^{(j)}) +\bfvarepsilon^{(j)}$, for $j = 1,\ldots,J$
        \State \hspace*{2ex} obtain $\bflambda_{\rm opt}^{(j)}$ (e.g., solve ~\eqref{eq:lambdatrue})
        \State \hspace*{2ex} set up DNN $\widehat \bfPhi$
        \State \hspace*{2ex} use training data $\left\{\, \bfb^{(j)},\bflambda_{\rm opt}^{(j)}\,  \right\}_{j= 1}^J$ to compute network parameters $\widehat \bftheta$ as in \eqref{eq:saa}
        \vspace*{2ex}
        \State \emph{online phase}
        \State \hspace*{2ex} obtain new data $\bfb$
        \State \hspace*{2ex} propagate $\bfb$ through the learned network to get $\widehat \bflambda = \widehat \bfPhi(\bfb; \ \widehat \bftheta)$
        \State \hspace*{2ex} compute inverse solution $\widehat\bfx (\widehat \bflambda)$ in~\eqref{eq:loss}
    \end{algorithmic}
    \caption{Learning regularization parameters via DNNs}\label{alg:nn}
\end{algorithm}

\begin{figure}[bthp]
\begin{center}
\includegraphics[width=1.0\textwidth]{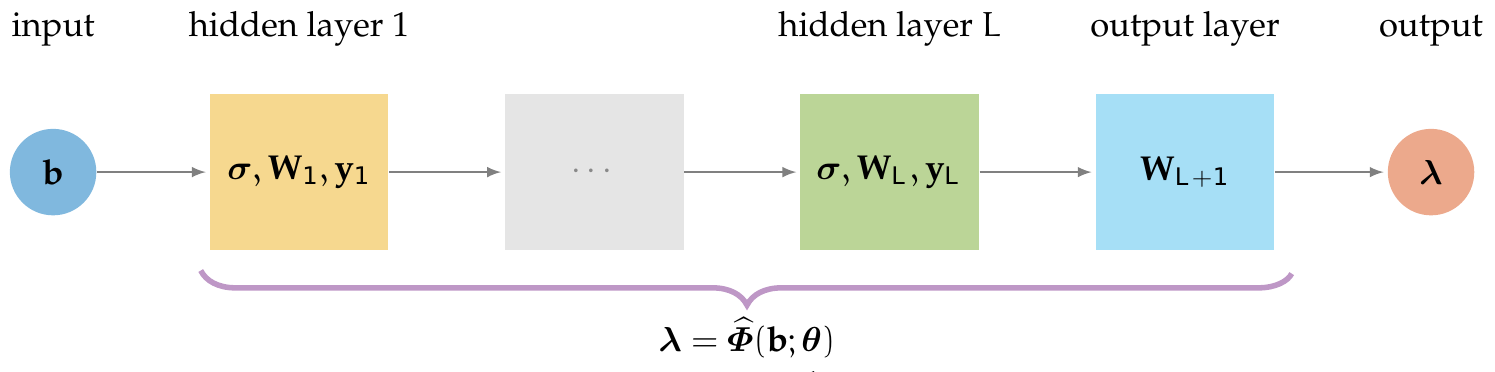}
\end{center}
\caption{Illustration of a DNN $\widehat\bfPhi(\bfb; \ \bftheta)$ with $L$ hidden layers. An input $\bfb$ is mapped by the network $\widehat\bfPhi$ onto an output $\bflambda$ given weights $\bfW_\ell$ and biases $\bfy_\ell$ for each layer. All weight and biases terms constitute the network parameter $\bftheta$.}
\label{fig:nn}
\end{figure}

In essence, our approach constructs a surrogate model using feedforward DNNs. Surrogate modeling methods are popular techniques used in scientific computing, where an approximate, trained model replaces the original model.  The surrogate model can be used for predicting outputs in unexplored situations or for reducing overall computational complexity \cite{gramacy2020surrogates}.
In feedforward models, information flows in one direction through a series of hidden layers of computation to produce the output. The connection between these hidden layers form a network, which is often represented as a directed graph. An illustration of a DNN is given in Figure~\ref{fig:nn}. Here, for simplicity of presentation, we illustrate and discuss fully connected neural networks.

Let's assume there exists a continuous \emph{target function} $\bfPhi:\bbR^{m} \to \bbR^{p}$, mapping observations $\bfb \in \bbR^{m}$ onto the regularization parameters $\bflambda\in \bbR^{p}$.  Our goal --- in a supervised machine learning approach ---  is to find a neural network $\widehat \bfPhi$ approximating the target function $\bfPhi$.  We define a \emph{fully-connected feedforward neural network} as a parameterized mapping $\widehat \bfPhi:\bbR^{m} \times \bbR^{q} \to \bbR^{p}$ with
\begin{equation} \label{eq:2.20}
	 \widehat \bfPhi(\bfb;\ \bftheta) = \bfvarphi_{L+1}(\bftheta_{L+1}) \circ \dots \circ \bfvarphi_1 (\bftheta_{1}) (\bfb),
\end{equation}
where $\circ$ denotes the component-wise composition of functions $\bfvarphi_\ell:\bbR^{m_{\ell-1}} \times \bbR^{q_\ell} \to \bbR^{m_\ell}$ for $\ell = 1, \ldots, L+1$  ($m_0 = m$, and $m_{L+1} = p$). The vector $\bftheta = \begin{bmatrix} \bftheta_1\t, & \ldots, & \bftheta_{L+1}\t\end{bmatrix}\t \in\bbR^q$ is a composition of layer specific parameters $\bftheta_\ell$ defining so-called \emph{weights} $\bfW_\ell$ and \emph{biases} $\bfy_\ell$, i.e., $\bftheta_\ell = \begin{bmatrix} \vec{\bfW_\ell}\t & \bfb_\ell\t \end{bmatrix}\t$, where $\bfW_\ell \in\bbR^{m_
\ell \times m_{\ell-1}}$ and $\bfb_\ell \in \bbR^{m_\ell} $. The functions $\bfvarphi_\ell$ are given by
\begin{equation}
    \bfvarphi_\ell(\bftheta_\ell)(\bfb_{\ell-1}) = \bfsigma_\ell( \bfW_\ell \bfb_{\ell-1} + \bfy_\ell ),
\end{equation}
where $\bfsigma_\ell:\bbR^{m_\ell} \to \bbR^{m_\ell}$ are typically nonlinear \emph{activation functions}, mapping inputs arguments point-wise onto outputs with limiting range.  Note that for the \emph{output layer} $\bfvarphi_{L+1}(\bftheta_{L+1})$, we assume a linear transformation with no bias term, $\bfvarphi_{L+1}(\bftheta_{L+1})(\bfb_L) = \bfW_{L+1} \bfb_{L}$.

The architecture or design of the neural network $\widehat \bfPhi$ is determined by the choice in the number of hidden layers $L$ \footnote{A neural network is considered deep if the network exceeds three layers including the input and output layer.}, the width of the layers (the size of $m_\ell$), and the type of the activation functions $\bfsigma_\ell$. For instance, if all $\bfvarphi_\ell$ are chosen to be the identity, then $\widehat \bfPhi$ becomes a linear transformation. Such networks provide excellent efficiency and reliability, but are not generally suited to approximate nonlinear mappings $\bfPhi$, see \cite{goodfellow2016machine}. Other choices for $\bfvarphi_\ell$ include differentiable functions, such as the logistic sigmoid and hyperbolic tangent function. While suitable for approximating smooth functions, computational inefficiencies during the training process work to their disadvantage \cite{goodfellow2016machine}. A popular choice is the \emph{rectified linear unit}, i.e., ${\rm ReLU}(x) = \max(0,x)$, hence $\bfsigma_\ell(\bfx) = \begin{bmatrix}{\rm ReLU}(x_1), & \ldots, & {\rm ReLU}(x_{m_\ell}) \end{bmatrix}\t$. The computationally efficient ${\rm ReLU}$ function is commonly used despite its non-differentiability. Practically, it has been observed that gradient based training methods are not impacted. The large number of degrees of freedom when defining the architecture of neural networks provides versatility and flexibility in approximating different types of functions. Approximation quality of the network $\widehat \bfPhi$ is application dependent and depends on the properties and complexity of the underlying function $\bfPhi$. In this regard, universal approximation properties for neural networks have been established, see for instance \cite{hornik1989multilayer,hornik1990universal,cybenko1989approximation}.

For imaging applications, it is appropriate for the DNN to integrate two dimensional distance structures into the design of the neural network.  In Section \ref{sec:numerics} we consider 2D convolutional neural networks, see~\cite{goodfellow2016machine}, where input data is convolved by kernels or filters. The weights of fully connected layers $\bfW_\ell$ become low dimensional filter factors $\bfW_\ell$ (and bias terms $\bfy_\ell$), which have the advantages of integrating the 2D structure of the problem and of reducing the ill-posedness of the learning of the network, thus allowing for deeper networks.

Determining a parameter set $\widehat\bftheta$ that leads to an accurate approximation $\widehat \bfPhi(\bfb; \widehat\bftheta) \approx \bfPhi(\bfb)$ is the key element of supervised learning via DNNs (see line 6 in Algorithm \ref{alg:nn}). We assume that training data $\left\{ \bfb^{(j)}, \bflambda_{\rm opt}^{(j)} \right\}_{j = 1}^J$ comprising of inputs $\bfb^{(j)}\in \bbR^m$ and corresponding outputs $\bflambda_{\rm opt}^{(j)}\in \bbR^p$ are available, and we select a cost function indicating the performance of $\widehat \bfPhi$. For \emph{regression} type problems, a common choice is the mean squared loss function $D(\cdot) = \norm[2]{\cdot}^2$.  Then the goal is to solve \eqref{eq:saa} to obtain the learned network parameters $\widehat \bftheta$.  However there are two considerations.  First, to prevent overfitting towards the training data, it is common to include an additional regularization term, e.g., $\calL(\bftheta) = \alpha^2||\bftheta||_2^2$.  For instance, we can solve a regularized problem,
\begin{equation}\label{eq:dnnlearn}
    \widehat \bftheta =\argmin_\bftheta \ \frac{1}{2J} \sum_{j=1}^J \norm[2]{\widehat\bfPhi(\bfb^{(j)}; \bftheta) - \bflambda_{\rm opt}^{(j)} }^2 + \calL(\bftheta).
\end{equation}
Indeed, the learning problem~\eqref{eq:dnnlearn} is itself a nonlinear inverse problem \cite{goodfellow2016machine}.   Second, there are various optimization methods that can be used to solve~\eqref{eq:dnnlearn}.  An intense amount of research in recent years has focused on the development of efficient and effective solvers for solving optimization problems like \eqref{eq:dnnlearn}.  \emph{Stochastic approximation} (SA) methods are iterative minimization approaches where a small subset of samples from the training set (e.g., a randomly chosen batch) is used at each iteration to approximate the gradient of the expected loss and to update the DNN weights \cite{shapiro2014lectures,robbins1951stochastic}. Common SA approaches include stochastic gradient descent and variants like ADAM \cite{kingma2014adam}. This approach is computational appealing for massively large datasets since only a batch of the data is needed in each step. However, slow convergence and the nontrivial task of selecting an appropriate step size or learning rate present major hurdles. As an alternative to SA methods, \emph{stochastic average approximation} (SAA) methods can be used, where a sample batch (or the entire training dataset) is used \cite{kleywegt2002sample}. One advantage is that deterministic optimization methods (e.g. inexact Newton schemes) can be used to solve the resulting optimization problem, but the main disadvantage is that a very large batch is typically required since the accuracy of the approximation improves with larger batch sizes.

\section{Learning the strength of regularization: one parameter}
\label{sec:oneparam}
Next, we investigate our proposed approach for the special but widely encountered problem where we seek \textit{one regularization parameter}, which represents the strength or amount of regularization. Without loss of generality, we consider a least squares loss function for the data fit. Let $\bflambda = \lambda \in \bbR$, then consider the regularized problem,
\begin{equation}\label{eq:1lambda}
  \widehat\bfx(\lambda) = \argmin_{\bfx} \ \norm[2]{\bfA(\bfx) -\bfb}^2 +\lambda^2 \calR(\bfx),
\end{equation}
where $\calR(\bfx)$ is a regularization functional that only depends on $\bfx$. In this case, the value of the regularization parameter $\lambda$ determines the weight or strength of the regularization term.  Another interpretation of $\lambda$ is that it represents the noise-to-signal ratio, which can be derived from a Bayesian perspective, see e.g., \cite{calvetti2007introduction,bardsley2018computational}.
In this section, we begin with an investigation on the use of a neural network to approximate the mapping from observation vector $\bfb$ to optimal regularization parameter $\lambda$ for the standard Tikhonov case.  Then, we address more general regularization terms and approaches.

\subsection{Standard Tikhonov regularization}
\label{sub:Tikhonov}
The standard Tikhonov problem, also referred to as ridge regression, is often used to solve linear inverse problems, see \cite{bardsley2018computational,hansen2010discrete,tenorio2017introduction}, where the regularized solution has the form,
\begin{equation}\label{eq:tik}
  \widehat\bfx(\lambda) = \argmin_{\bfx} \ \norm[2]{\bfA\bfx -\bfb}^2 +\lambda^2 \norm[2]{\bfx}^2,
\end{equation}
with $\bfA \in \bbR^{m \times n}$. An adequate choice of the regularization parameter $\lambda$ is paramount, since a parameter that is too small may lead to erroneous solutions and a parameter that is too large may lead to overly smoothed solutions.  Assuming $\bfx_{\rm true}$ is known, an optimal (balanced) regularization parameter can be computed as in~\eqref{eq:lambdatrue} with $\widehat \bfx(\lambda)$ given by~\eqref{eq:tik}.
In practice, $\bfx_{\rm true}$ is not available, and there are various regularization parameter selection methods to estimate $\lambda_{\rm opt}$.
Prominent methods include the discrepancy principle (DP), the generalized cross-validation (GCV) method, the unbiased predictive risk estimator (UPRE), and the residual periodogram, see e.g., \cite{hansen2010discrete,bardsley2018computational} and references therein.  If an estimate of the noise variance $\sigma^2$ is available, two popular approaches include DP and UPRE.  The DP parameter is computed by solving the root finding problem,
\begin{equation}
\label{eq:dp}
    \lambda_{\rm dp}  =  \left\{ \lambda : \norm[2]{\bfA\widehat\bfx(\lambda)  -\bfb}^2 - m \sigma^2 =0  \right\},
\end{equation}
and the UPRE parameter
can be computed as,
\begin{equation}\label{eq:upre}
    \lambda_{\rm upre} = \argmin_{\lambda} \ U(\lambda) =  \norm[2]{\bfA\widehat\bfx(\lambda) -\bfb}^2 + 2\sigma^2 \trace{\bfA\bfZ(\lambda)}.
\end{equation}
Here, $\trace{\bfB}$ denotes the trace of a matrix $\bfB$ and $\bfZ(\lambda) = (\bfA\t\bfA + \lambda^2 \bfI_n)^{-1} \bfA\t$. A common parameter choice method that does not require an estimate of the noise variance is the GCV method, which is based on leave-one-out cross validation \cite{bardsley2018computational}.  The goal is to find a regularization parameter that solves the following optimization problem,
\begin{equation}\label{eq:gcv}
\lambda_{\rm gcv}  = \argmin_{\lambda} \ G(\lambda) = \frac{m\norm[2]{\bfA\widehat\bfx(\lambda) -\bfb}}{(\trace{\bfI_m-\bfA\bfZ(\lambda)})^2}.
\end{equation}
For each of the described methods, we assume a unique solution exists.
Many of these methods have sound theoretical properties (e.g., statistical derivations) and can lead to favorable estimates (e.g., providing unbiased estimates of $\lambda$). However, each approach has some disadvantages.
For example, DP and UPRE require estimation of the noise variance $\sigma^2$ \cite{donoho1995noising}, and the computational costs to minimize the UPRE and GCV functional are significant, since computing $U(\lambda)$ and $G(\lambda)$ may involve $\calO(n^3)$ floating point operations.  For small problems or for problems where the singular value decomposition (SVD) of $\bfA$ is available, the SVD can be used to significantly reduce the costs to compute parameters \eqref{eq:dp}, \eqref{eq:upre}, and \eqref{eq:gcv}.  Furthermore, in recent years, various methods have aimed to reduce computational costs, e.g., through trace estimation and other randomized linear algebra approaches \cite{luiken2020comparing}. However, for many large-scale problems, the burden of computing a suitable regularization parameter $\lambda$ still remains.  Indeed, the computational cost to compute an estimate of $\lambda$ using standard techniques oftentimes far outweighs the cost of solving the inverse problem~\eqref{eq:tik} itself.
One remedy is to consider hybrid projection methods that combine an iterative projection method with a variational regularizer so that the regularization parameter can be automatically tuned on a much smaller, projected problem \cite{o1981bidiagonalization,bjorck1988bidiagonalization,chung2008weighted}.  Nevertheless, selecting an appropriate regularization parameter using existing approaches may still be difficult or costly in practice.  In the following, we consider the use of DNN-predicted regularization parameters and begin with a small test problem to provide comparisons to existing parameter choice methods.

\begin{figure}[t]
    \begin{center}
    \includegraphics[width=0.45\textwidth]{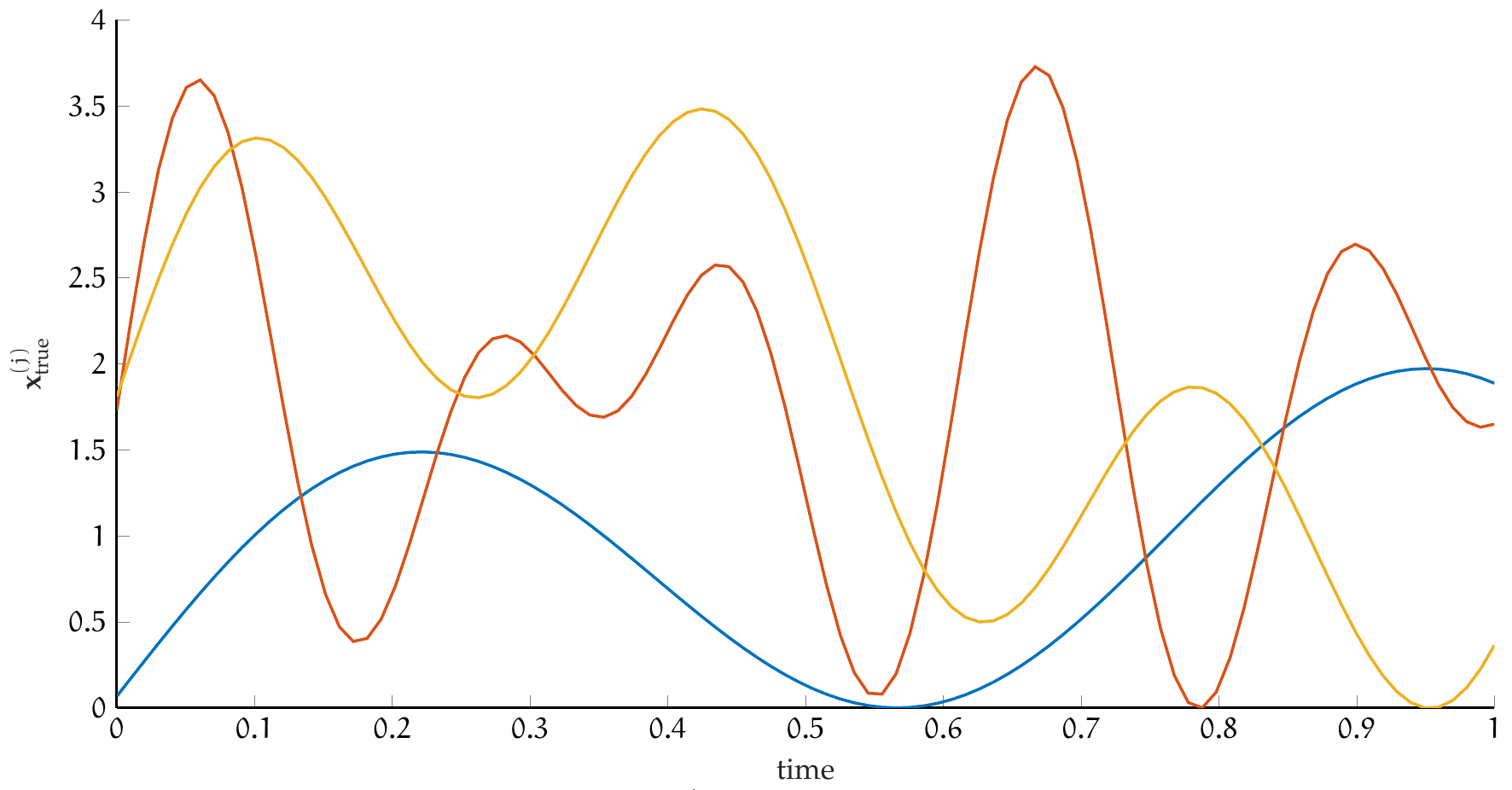}
    \includegraphics[width=0.45\textwidth]{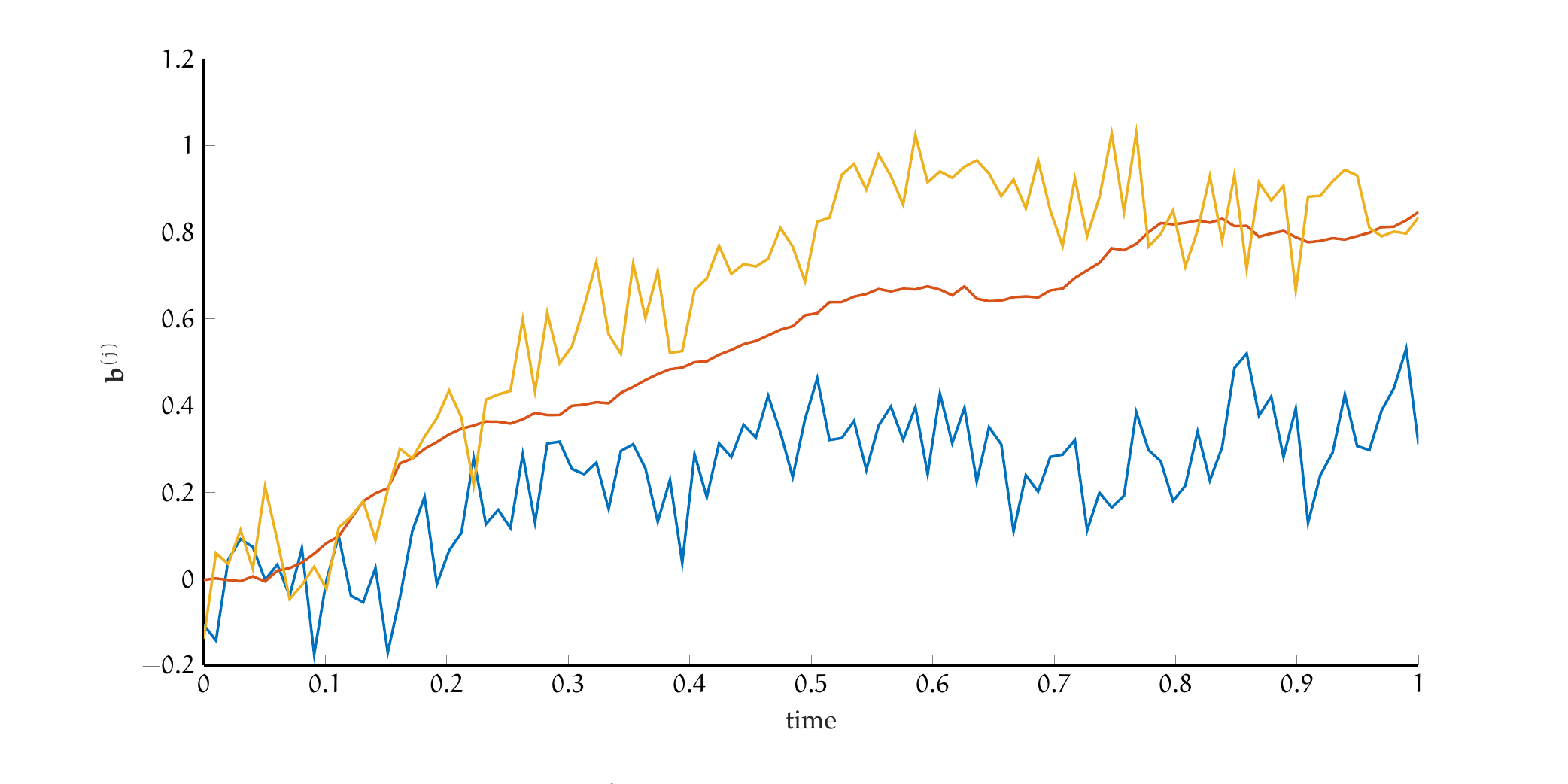}
    \caption{For the inverse heat condition problem, the left panel shows three sample temperature signals $\bfx_{\rm true}$ at location $0$ in time. The right panel depicts the corresponding noisy temperature observations at location $1$, with different noise characteristics.}\label{fig:signal1}
    \end{center}
\end{figure}

Consider the \emph{classical inverse heat conduction problem}:  An unknown heat source is applied to the end of an insulated semi-infinite bar (at location $0$). Given noise contaminated temperature measurements $\bfb = [b(t_1),\ldots,b(t_m)]\t$ at time points $t_1,\ldots,t_m$ at location 1, the goal is to determine the temperature of the heat source $x(t)$ at any time $t$ at location $0$, see \cite{lamm2000survey}.  For the simulated problem, let $m=n=100$ and assume a heat source of the form $\bfx_{\rm true} = [x(t_1),\ldots,x(t_{100})]\t$ with $x(t) = \sin(2\pi r_1 t) + \sin(2\pi r_2t) +c$ at location $0$, where $r_1$ and $r_2$ are random parameters and $c$ is selected such that $x(t)\geq 0$. Let $\bfA\in\bbR^{n\times n}$ represent the forward operator, as computed in the {\tt regTools} Matlab toolbox~\cite{regTools}.
Then the synthetic data are generated as $\bfb = \bfA\bfx_{\rm true} + \bfvarepsilon$ with noise $\bfvarepsilon \sim \calN(\bfzero, \sigma^2\bfI_m)$ for some noise variance $\sigma^2$. By randomly selecting $\sigma^2$ to be uniformly distributed between $10^{-3}$ and $10^{-1}$ and randomly setting parameters $r_1$ and $r_2$, we generate a training set of size $J = 200,\!000$, see Figure~\ref{fig:signal1} depicting three samples. Following step 4 of Algorithm \ref{alg:nn}, we compute the optimal regularization parameter $\lambda_{\rm opt}^{(j)}$ for each of the samples.

We consider two network designs: a deep network and a shallow network.
\begin{description}
    \item[DNN] For the DNN, we utilize a fully connected feedforward network $\widehat \bfPhi(\mdot;\bftheta):\bbR^{100} \to \bbR$ with five hidden layers of widths 75, 50, 25, 12, and 6; hence, the network parameters are contained in $\bftheta \in \bbR^{13,046}$.
    Using the training data $\left\{\bfb^{(j)} ,\lambda_{\rm opt}^{(j)}\right\}_{j = 1}^J$, we compute an estimate of the DNN network parameters denoted as $\widehat \bftheta$ by solving \eqref{eq:saa}
    (c.f., step 6 of Algorithm \ref{alg:nn}).
    Since this is a fairly small problem, the optimization of the regression problem~\eqref{eq:dnnlearn} is performed on the entire data set (i.e., an SAA approach) using a Levenberg-Marquardt method with regularization functional $\calL \equiv 0$.
    \item[ELM] As a second neural network, we consider an extreme learning machine (ELM) which consists only of an output layer containing weights $\bfw\in\bbR^{100}$ and a bias term, $y\in \bbR$, see \cite{huang2006extreme}. The ELM model corresponds to a simple assumption that there is a linear relationship between the input data $\bfb$ and the output $\lambda$, i.e., $\lambda = \bfw\t\bfb+y$. With regression, training this shallow neural network (i.e., estimating $\bfw$ and $y$) simplifies to a linear least squares problem which can be solved efficiently using an iterative method (e.e., {\tt lsqr} \cite{paige1982lsqr}). Let $\widehat\bfw,\widehat y$ denote the computed ELM parameters.
\end{description}

\begin{figure}[t]
    \centering
    \includegraphics[width=0.45\textwidth]{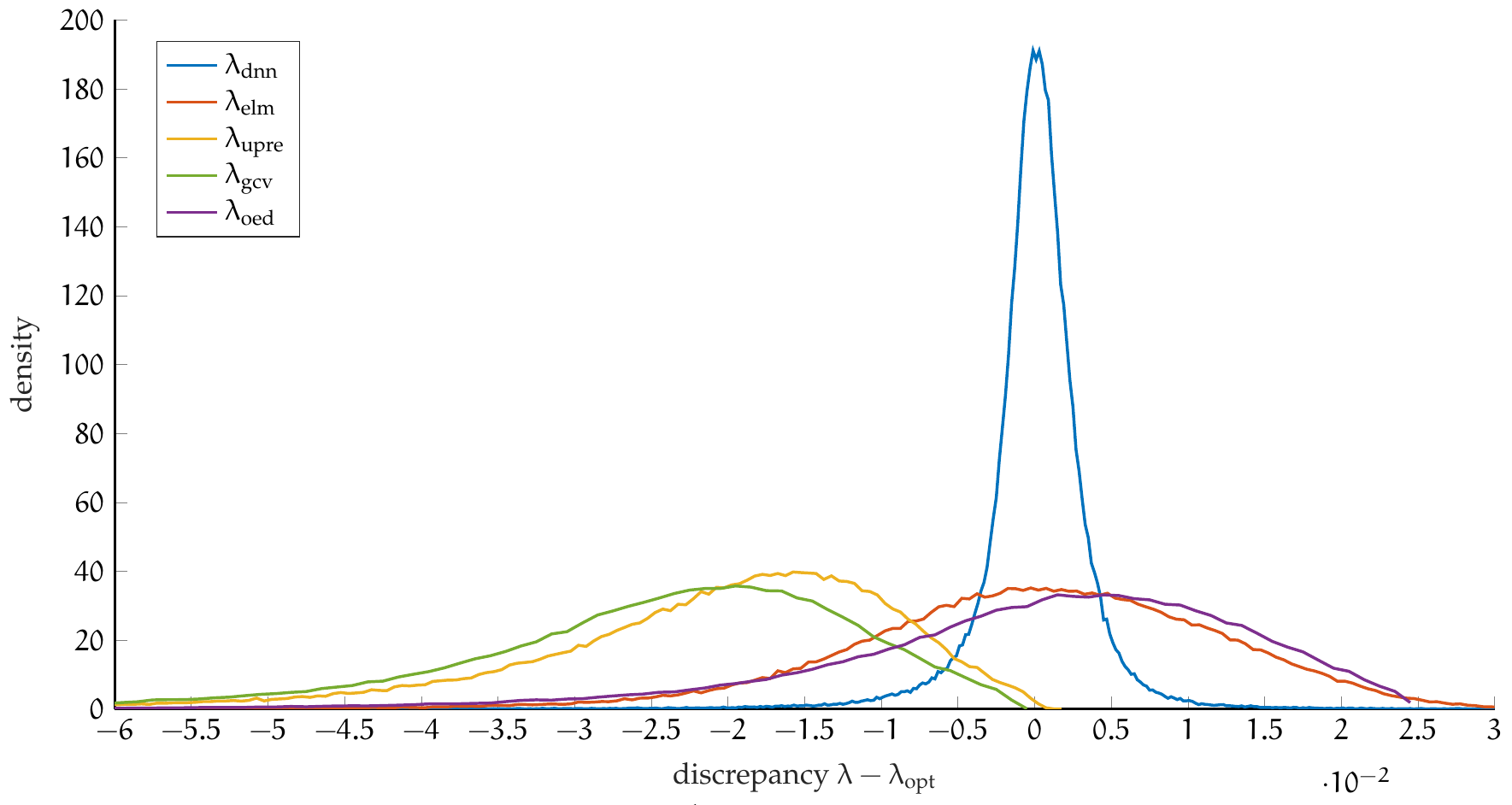}
    \includegraphics[width=0.45\textwidth]{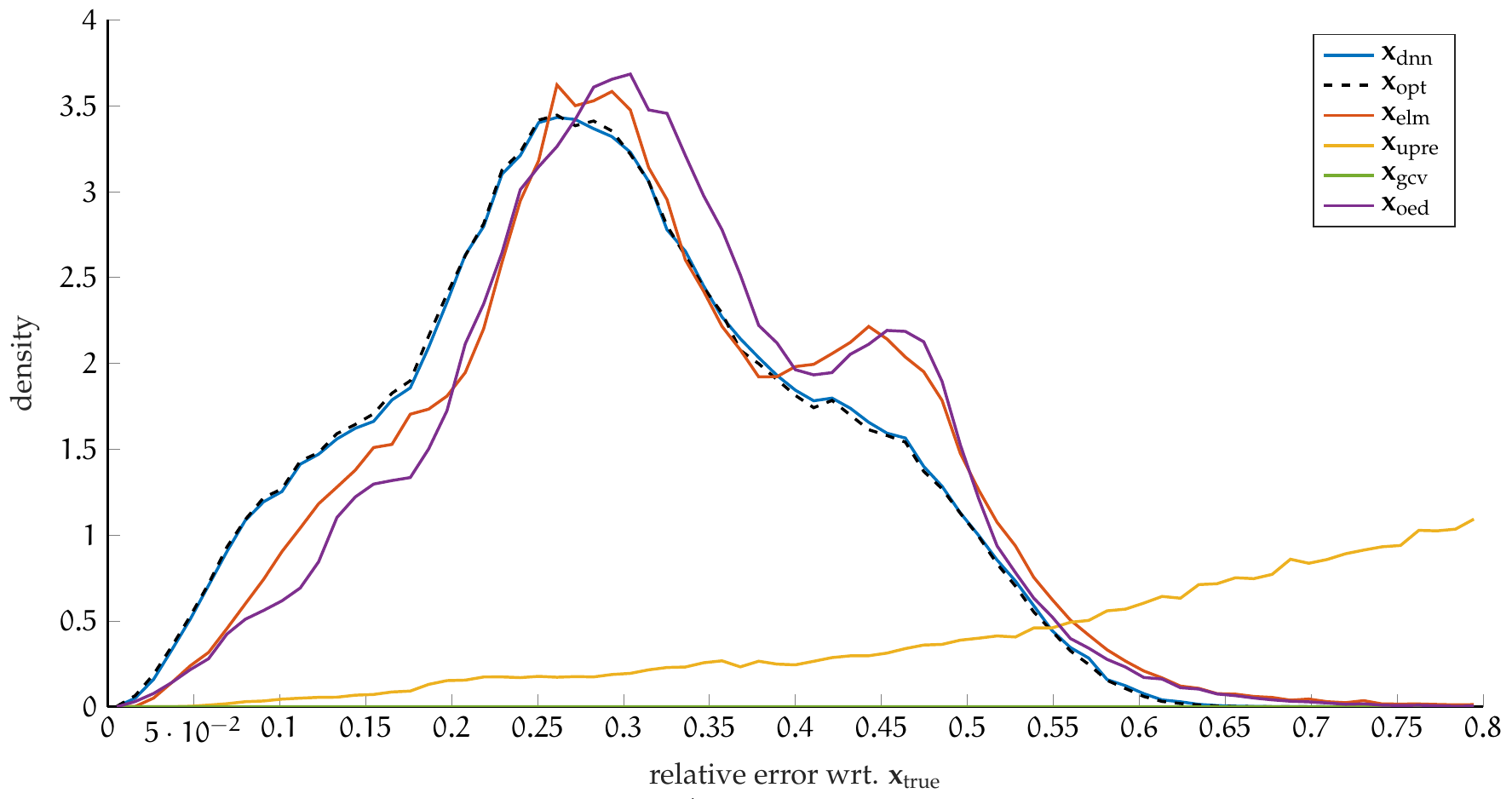}
    \caption{The left panel displays the distribution of the discrepancy between the estimated regularization parameter $\lambda$ and $\lambda_{\rm opt}$ for various parameter choice methods for 200,000 validation data. On the right we provide the distribution of relative reconstruction error norms for various parameter choice methods.}\label{fig:regpar}
\end{figure}

Next we describe the online phase.  Using the same randomized procedure described above, we generate a validation set of size $J = 200,\!000$.
For each sample from the validation set, $\bfb^{(j)}$, we compute the DNN predicted regularization parameter as $\lambda^{(j)}_{\rm dnn} = \bfPhi(\bfb^{(j)};\widehat\bftheta)$ and the ELM predicted regularization parameter as $\lambda^{(j)}_{\rm elm} = \widehat\bfw\t \bfb^{(j)}+\widehat y$. The corresponding DNN and ELM reconstructions are denoted as $\bfx^{(j)}_{\rm dnn} = \widehat \bfx(\lambda^{(j)}_{\rm dnn})$ and $\bfx^{(j)}_{\rm elm}=\widehat \bfx(\lambda^{(j)}_{\rm elm})$ respectively. To evaluate the performance of the network computed parameters, for each sample from the validation set, we compute the discrepancy to the optimal regularization parameter $\lambda_{\rm opt}^{(j)}$ and the relative error norm of the reconstruction $\widehat\bfx$ with respect to $\bfx_{\rm true}$, $\norm[2]{\widehat\bfx - \bfx_{\rm true}}/\norm[2]{\bfx_{\rm true}}$.
In Figure~\ref{fig:regpar}, we provide distributions of the discrepancy in computed regularization parameter in the left panel and distributions of the relative errors in the temperature reconstructions in the right panel.  We observe that both of the network based approaches perform reasonably well.

For comparison, we also provide results corresponding to the UPRE-selected regularization parameter $\lambda^{(j)}_{\rm upre}$ and the GCV-selected regularization parameter $\lambda^{(j)}_{\rm gcv}$ (see equations~\eqref{eq:upre} and~\eqref{eq:gcv}). We exclude results for DP due to significant under-performance.  We also provide a comparison to an OED approach where regularization parameter $\lambda_{\rm oed}$ is computed by minimizing~\eqref{eq:OED} for the training set, and that \textit{one} parameter is used to obtain all reconstructions $\bfx^{(j)}_{\rm oed} = \widehat \bfx(\lambda_{\rm oed})$ for the validation set.

Recall $\lambda_{\rm opt}^{(j)}$ corresponds to the theoretical optimal performance which cannot be obtained in real world problems. From Figure~\ref{fig:regpar} we observe that both GCV and UPRE are under-performing and underestimate the optimal regularization parameter. Compared to standard methods, the OED approach performs significantly better at estimating $\lambda_{\rm opt}^{(j)}$.  ELM generates slightly better results then the OED approach. However, results from the DNN are virtually indistinguishable from results obtained using the optimal regularization parameters $\lambda_{\rm opt}^{(j)}$ and therefore perform extremely well.

In summary, both network predicted approaches (DNN and ELM) outperform existing methods. While ELMs have slight computational advantages compared to DNNs, we will see in Section~\ref{sec:numerics} that DNNs can better predict regularization parameters, thereby resulting in smaller reconstruction errors than ELMs. Notice that the network predicted parameters are very close to the optimal ones, which corresponds to very accurate regularized solutions. Furthermore, they are very cheap to compute.  That is, given new data, the network-predicted regularization parameter can be computed with one feedforward evaluation of the neural network. In fact, by shifting the main computational costs, i.e., training the neural network, to the offline phase, the computational complexity of the online phase is significantly reduced. Notice also that neural networks are more versatile than OED approaches, since in the OED approach only one regularization parameter is computed and that parameter is strongly dependent on the design choice as well as the training data.  For a squared loss design function, the computed parameter is only good \emph{on average} for the training data of a given problem~\cite{Pukelsheim2006,Atkinson2007}.  Thus, a major drawback of OED methods is limitations in generalizability (e.g., to observations with different noise levels or other features). On the contrary, network learned parameters obtained in an online phase are tailored to the new data.

\subsection{General regularization}
The approach to learn a regularization parameter for Tikhonov regularization via training of neural networks described in Section \ref{sub:Tikhonov} can be extended to more general regularization terms and approaches. Indeed, the problem of estimating a good regularization parameter for \eqref{eq:1lambda} becomes significantly more challenging for nonlinear problems and for non-traditional, nonlinear regularization terms. A common approach is to spend a significant effort to compute a good regularization parameter a priori and then to solve the optimization problem for fixed regularization parameter \cite{engl1996regularization}.  This can be very expensive, requiring many solves for different parameter choices in the online phase.  For nonlinear inverse problems, another approach uses a two-stage method that first reduces the misfit to some target misfit value and second to keep the misfit fixed and to reduce the regularization term. Although very popular in practice, this approach is not guaranteed to converge (in fact diverging in some cases) and appropriate safety steps and ad hoc parameters are needed \cite{parker1994geophysical,constable1987occam}.  For nonlinear least-squares problems with a Tikhonov term, Haber and Oldenburg in \cite{haber2000gcv} combine a damped Gauss-Newton method for local regularization with a generalized cross-validation method for selecting the global regularization parameter, but the overall scheme can still be costly.

For more general (non-Tikhonov) regularizers, selecting a regularization parameter can be computationally costly even for linear problems \cite{liao2009selection,lin2010upre,wen2011parameter,langer2017automated}. Total variation (TV) regularization~\cite{borsic2009vivo,rudin1992nonlinear} is a common approach, where the penalty term or regularizer takes the form
\begin{equation}
\label{eq:TV}
    \calR(\bfx,\lambda) = \lambda TV(\bfx)
\end{equation}
with regularization parameter $\lambda>0$ and $TV$ representing the total variation function. Anisotropic TV is often used when one seeks to promote sparsity in the derivative, i.e., partial smoothness.
Standard regularization parameter selection methods (e.g., DP, UPRE, and GCV approaches) are not easily extendable, although more elaborate methods based on these principles have been considered, e.g., \cite{wen2011parameter}. Sparsity-promoting regularizers based on $\ell_p$ regularization and inner-outer schemes for edge and or discontinuity preservation have gained popularity in recent years \cite{rudin1992nonlinear,arridge2019solving}, but selecting regularization parameters for these settings is not trivial.  Various extensions of hybrid projection methods to more general settings have been developed \cite{gazzola2014generalized, chung2019flexible,gazzola2020inner}. Such methods exploit iteratively reweighted approaches and flexible preconditioning of Krylov subspace methods in order to avoid expensive parameter tuning, but can still be costly if many problems must be solved. Furthermore, there are many works on supervised learning in an OED framework \cite{haber2003learning,haber2008numerical} for nonlinear problems \cite{horesh2010optimal} and general regularization terms (e.g., total variation \cite{de2017bilevel}, and sparsity \cite{haber2012numerical}).

Another common form of regularization is iterative regularization, where regularization is imposed by early termination of some iterative (often projection based) approach, applied to the model-data misfit term of~\eqref{eq:loss}. Iterative regularization methods are widely used for solving large-scale inverse problems, especially nonlinear ones with underlying partial differential equations (PDEs) (e.g., parameter identification in electrical impedance tomography), due to their ability to handle more complex forward models.  For example, most iterative methods only require the operation of the forward model $\bfA(\bfx)$ at each iteration, which is ideal for problems where the forward models can only be accessed via function evaluations.  For linear problems, matrices $\bfA$ and $\bfA\t$ may be too large to be constructed, but evaluations can be done cheaply (e.g., by exploiting sparsity or structure). The main challenge with iterative regularization is determining a good stopping iteration, which is complicated due to a phenomenon called \textit{semi-convergence}. Many iterative Krylov subspace methods when applied to inverse problems exhibit semi-convergence behavior, where during the early iterations the solution converges to the true solution, but at some point, amplification of the noise components in the approximate solution lead to divergence from $\bfx_\mathrm{true}$ and convergence to the corrupted and undesirable naïve solution.  This change occurs when the Krylov subspace begins to approximate left singular vectors corresponding to the small singular values.  For a simulated problem where we know $\bfx_\mathrm{true}$, one way to visualize semi-convergence is to plot the relative reconstruction error norms per iteration, which exhibits a ``U''-shaped plot; see the black line in Figure~\ref{fig:relerrors}.  Stopping the iterative process too early can result in images that are too smooth, and stopping too late can result in severely degraded reconstructions. The stopping iteration plays the role of the regularization parameter, and it can be very challenging to determine appropriate stopping criteria.

By defining a neural network that maps the observed data $\bfb$ to an optimal regularization parameter or an optimal stopping iteration, Algorithm \ref{alg:nn} can be used to efficiently estimate a parameter defining the strength of regularization for different regularization approaches. In the next section, we provide various numerical experiments from image processing to show that DNNs are well suited to approximate the strength of regularization, as well as other parameters describing regularity.

\section{Numerical Experiments}
\label{sec:numerics}
In this section, we provide several numerical examples from image processing that demonstrate the performance of our proposed approach for learning regularization parameters.  In Section~\ref{sec:tomo}, we consider a tomographic reconstruction example where a DNN is used to approximate the mapping from observation to the regularization parameter for TV regularization.  In Section \ref{sec:deblur}, we consider an image deblurring example where the goal is to detect outlines of inclusions in density fields from blurred images.  We demonstrate how DNNs can be used to approximate both the TV regularization parameter and the parameter quantifying the degree of an object's regularity.  A third example is provided in Section~\ref{sec:iter}, where we consider an inverse diffusion problem where regularization is enforced by early stopping of an iterative method (i.e., iterative regularization), and we train a DNN to estimate an appropriate stopping iteration.

We remark that for all of our experiments, we noticed that the network learning process was very robust in regards to the number of layers, the width of the layers, and the overall design of the network.  We attribute this to the fact that we have only a few output parameters.

\subsection{Computerized Tomography Reconstruction} \label{sec:tomo}
\begin{figure}[bthp]
	\centering
    \includegraphics[width=0.5\textwidth]{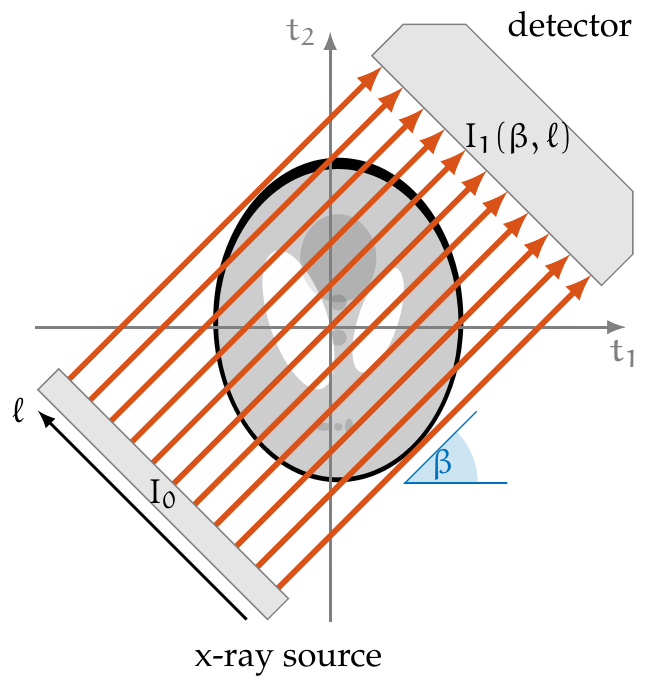}
		\caption{Illustration of 2D x-ray parallel-beam tomography setup, modified from \cite{ruthotto2018optimal}.}
		\label{fig:tomo:tomo2d}
\end{figure}

Computerized tomography (CT) is a widely used imaging technique for imaging sections or cross-sections of an object using penetrating waves.  For example, in biomedical imaging, x-rays are passed through some medium and, dependent on the properties of the material, some of the energy from the x-rays is absorbed.  Detectors with multiple bins measure the intensity of the x-rays emitted from the source (i.e., energy from x-rays that pass through the medium). See Figure~\ref{fig:tomo:tomo2d} for a visualization in 2D.

By rotating the x-ray source and/or the detectors around the object, measurements are collected from different angles.  These measurements are contained in the sinogram. CT reconstruction is a classic example of an \emph{inverse problem}, where the goal is to \emph{infer} the energy absorbency of the medium from the sinogram.

With appropriate simplifications, the discrete CT problem can be stated as \eqref{eq:invprob} with $\bfA(\bfx_{\rm true}) = \bfA\bfx_{\rm true}$, see \cite{natterer2001mathematics,kak2002principles}.
In this example, $\bfx_{\rm true} \in \bbR^{16,384}$ represents a discretized version of the observed medium ($128 \times 128$), the forward operator $\bfA$ represents the radon transform, $\bfb$ represent the observed intensity loss between detector and x-ray source (with $181$ parallel rays over $180$ equidistant angles $\theta$), and $\bfvarepsilon$ reflects noise in the data. Tomography problems are ill-posed inverse problems and regularization is required. We consider TV regularization \eqref{eq:TV}.

\begin{figure}[bthp]
    \centering
    \begin{tabular}{ccc}
        \includegraphics[width= 0.20\textwidth]{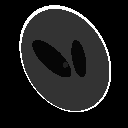} &
        \includegraphics[width= 0.20\textwidth]{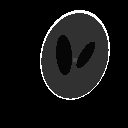} &
        \includegraphics[width= 0.20\textwidth]{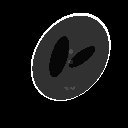}\\[1.5ex]
        \includegraphics[width= 0.20\textwidth]{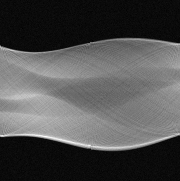} &
        \includegraphics[width= 0.20\textwidth]{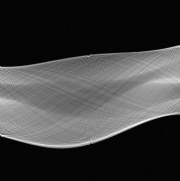} &
        \includegraphics[width= 0.20\textwidth]{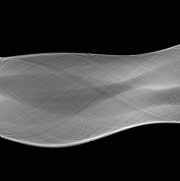}\\[1.5ex]
        \includegraphics[width= 0.20\textwidth]{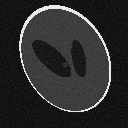} &
        \includegraphics[width= 0.20\textwidth]{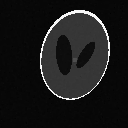} &
        \includegraphics[width= 0.20\textwidth]{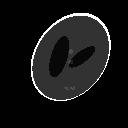}\\[1.5ex]
    \end{tabular}
    \caption{Samples from the tomography dataset: three training images $\bfx_{\rm true}^{(j)}$ (first row), noisy sinograms $\bfb^{(j)}$ (second row) with white noise level $2.58 \%$,   $1.37 \%$, and  $0.48 \%$, respectively, and optimal reconstructions $\bfx_{\rm opt}^{(j)}$ (last row).}
    \label{fig:tomodata}
\end{figure}

The goal of this experiment is to train a DNN to represent the mapping from sinogram to TV regularization parameter. First, we generate true images $\bfx_{\rm true}^{(j)}$ using the random Shepp-Logan phantom see~\cite{ruthotto2018optimal,randomSheppLogan}. Then, sinograms are computed as $\bfb^{(j)} = \bfA\bfx_{\rm true}^{(j)} +\bfvarepsilon^{(j)}$, where $\bfvarepsilon$ represents white noise with a noise level that is selected as uniform random between $0.1$\% and $5$\%. For example, a noise level of $5$\% corresponds to $\norm[2]{\bfvarepsilon}/\norm[2]{\bfA \bfx_{\rm true}}=0.05.$ The operator $\bfA$ is determined via parallel beam tomography \cite{hansen2018air}. For each of the sinograms, we solve~\eqref{eq:lambdatrue} using a golden section search algorithm to obtain $\lambda_{\rm opt}^{(j)}$.
This requires solving multiple TV regularization problems which is performed using the split-Bregman approach, see \cite{rudin1992nonlinear} for details.
We denote the corresponding optimal reconstruction by $\bfx_{\rm opt}^{(j)} = \widehat\bfx(\lambda_{\rm opt}^{(j)})$.
In Figure~\ref{fig:tomodata}, we provide three of the true images in the top row, the corresponding observed sinograms in the middle row, and the TV reconstructions corresponding to the optimal regularization parameter in the bottom row.

We select a training set of $24,\!000$ images and consider learning approaches to approximate the input-output mapping $\bfb^{(j)} \to \lambda_{\rm opt}^{(j)}$ for $j = 1,\ldots, 24,\!000$. Notice that the network inputs are image sinograms of size $181 \times 180$.  For this application we consider a convolutional neural network consisting of 4 single channel convolutional layers with $32\times 32$, $16\times 16$, $8\times8$, and $4\times4$ kernel weights, respectively and one bias term each. The convolutional layers are padded, and the kernel is applied to each pixel, i.e., the stride is set to~1. To reduce the dimensionality of the neural network we use an average pooling of $32\times 32$. We establish one fully connected layer with $4\times 22,\!350$ weights, plus $4$ bias terms, and a one dimensional output layer $1\times4$, plus one bias term. Each hidden layer has a ReLU activation function and with a regression loss there are $90,\!773$ parameters in $\bftheta$ defining the DNN $\widehat\bfPhi(\bfb;\bftheta)$. To estimate $\bftheta$ we utilize the ADAM optimizer with a learning rate of $10^{-4}$, while the batch size is set to $64$,~\cite{kingma2014adam}. We learn for $30$ epochs. Further, we also consider a shallow network and use an ELM design as described in Section~\ref{sec:oneparam}.

\begin{figure}[bthp]
    \begin{center}
    \includegraphics[width=0.7\textwidth]{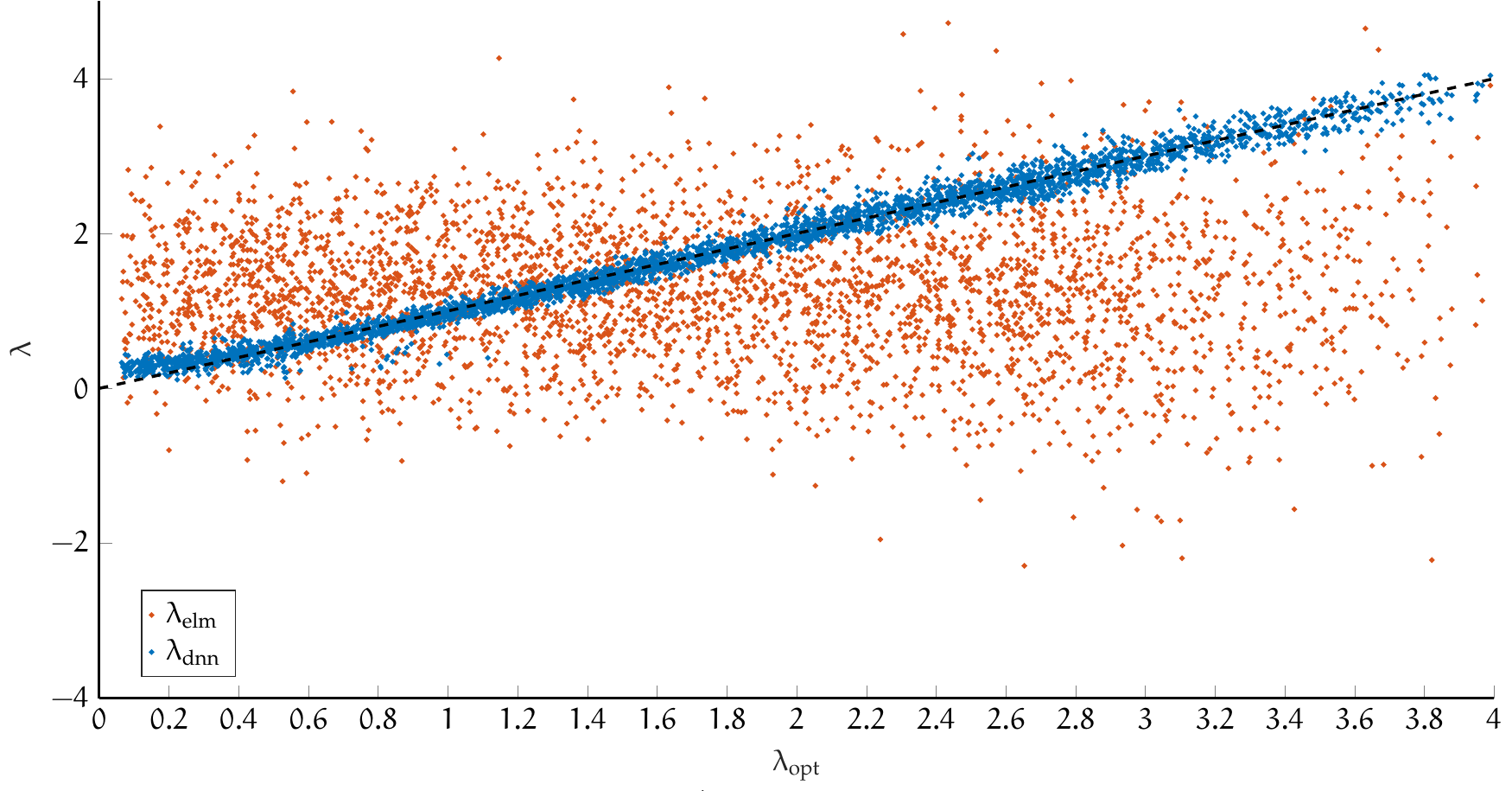}
    \caption{Scatter plot of network predicted regularization parameters $\lambda_{\rm dnn}$ and $\lambda_{\rm elm}$ versus the optimal regularization parameter $\lambda_{\rm opt}$ for the tomography reconstruction example. }\label{fig:scatterTomo}
    \end{center}
\end{figure}

A validation set of $3,\!600$ images are generated with the same properties as the training set. The scatter plot in Figure~\ref{fig:scatterTomo} illustrates the predictive performance of the DNN and the ELM networks on the validation set. For each sample from the validation data set, we compute the network predicted regularization parameters $\lambda_{\rm dnn}$ and $\lambda_{\rm elm}$ and plot these against the optimal regularization parameter $\lambda_{\rm opt}$.

\begin{figure}[bthp]
    \begin{center}
    \includegraphics[width=0.45\textwidth]{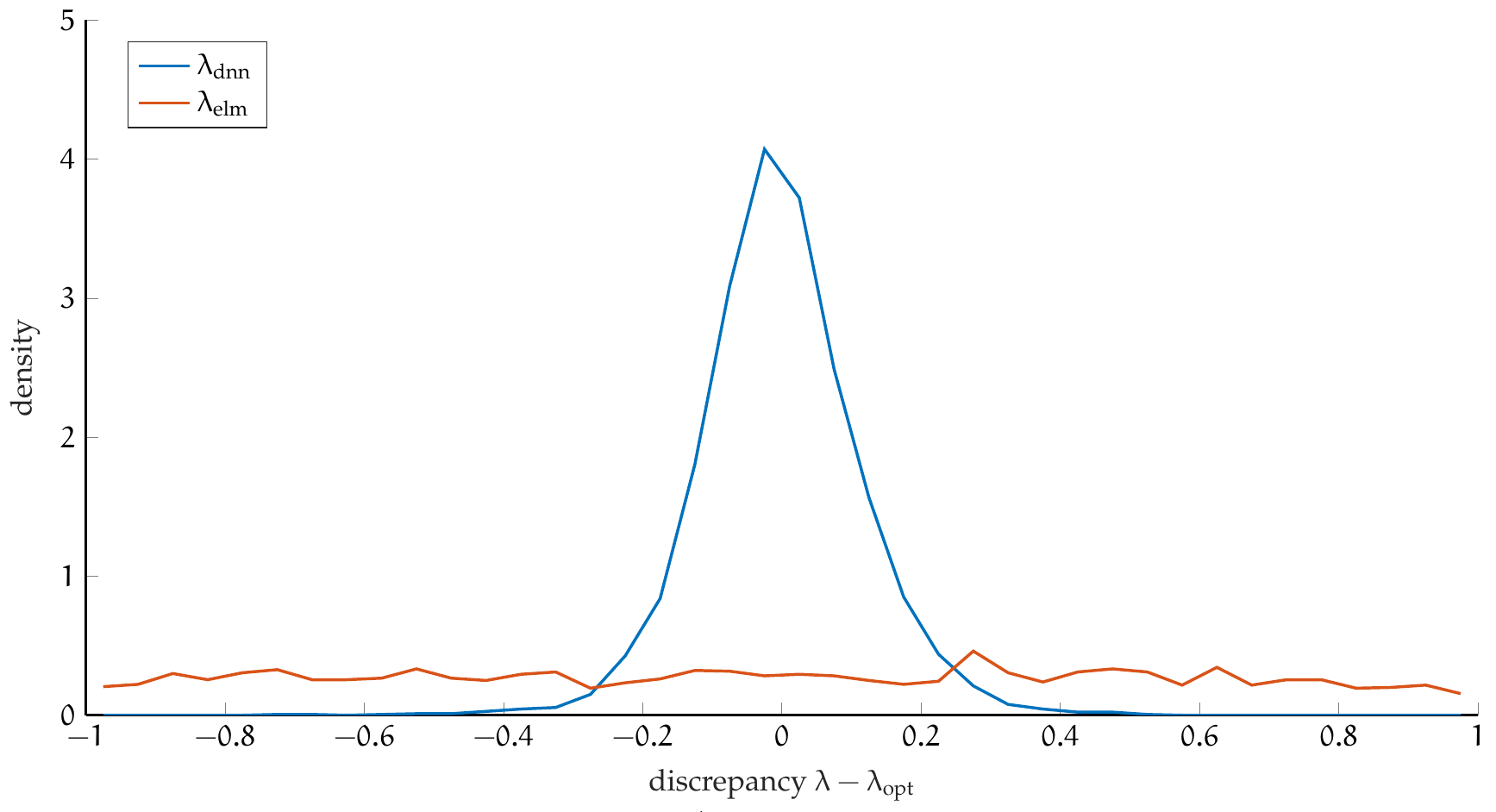}
    \includegraphics[width=0.45\textwidth]{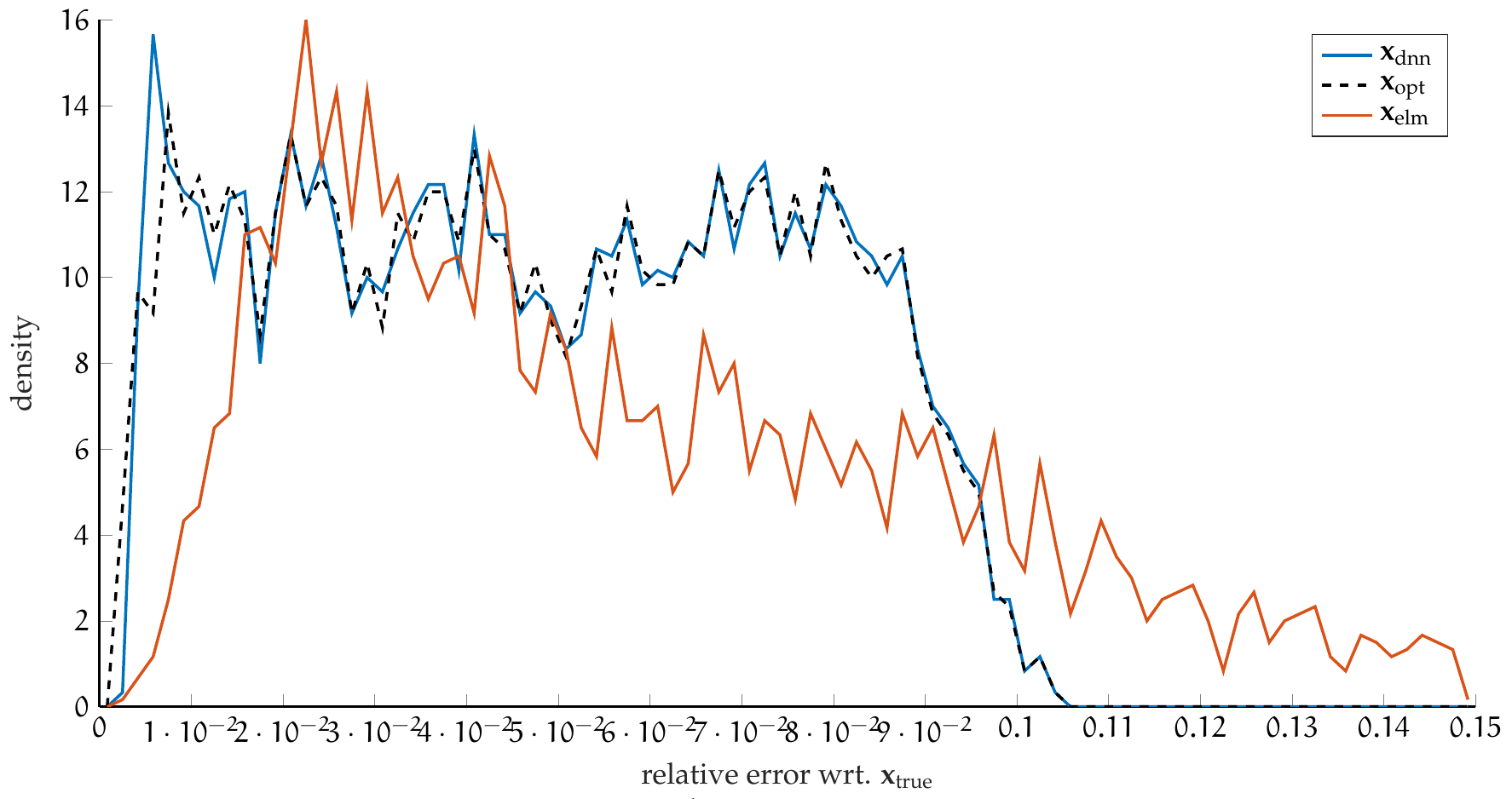}
    \caption{(left) Probability densities for the discrepancy between network predicted regularization parameters and the optimal regularization parameter for the tomography reconstruction example. (right) Probability densities for the relative reconstruction error norms.}\label{fig:resultsTomo}
    \end{center}
\end{figure}

While the scattered data of the ELM vary greatly, the scattered data of the DNN clusters around the identity line revealing the favorable predictive performance of the DNN.
Another way to visualize these results is to look at the discrepancy between the network predicted regularization parameters $\lambda_{\rm dnn}$ and $\lambda_{\rm elm}$ and the optimal regularization parameter.
The probability densities of these discrepancies are provided in the left panel of Figure~\ref{fig:resultsTomo} and further reveal the alignment between the DNN computed and the optimal regularization parameter.
Next we investigate the translation of the network predictions of the regularization parameters to the quality of the image reconstruction.  We compute relative reconstruction error norms as $\norm[2]{\widehat\bfx(\lambda) -\bfx_{\rm true}}/\norm[2]{\bfx_{\rm true}}$ for $\lambda_{\rm dnn}, \lambda_{\rm elm}$ and $\lambda_{\rm opt}$ and provide densities for the validation data in the right panel of Figure~\ref{fig:resultsTomo}. While the ELM predicted regularization parameter resulted in significant errors, partially due to large outliers, we observe that the DNN predicted regularization parameters resulted in near optimal TV reconstructions. In fact, we found that in a few instances the DNN predicted regularization parameter resulted in a smaller relative reconstruction error than the optimal, which revealed numerical errors in the computation of the optimal regularization parameter $\lambda_{\rm opt}$.

\subsection{Image deblurring with star-shaped inclusions}
\label{sec:deblur}

Digital imaging is an important tool in medicine, astronomy, biology, physics and many industrial applications and is frequently used to answer cutting edge scientific questions. Imperfections in imaging instruments and models can result in blurring and degradation of digital images. Post-processing methods for removing such artifacts from a digital image have been a topic of active research in the past few decades~\cite{Puetter2005,Pearson2018}. For the study in this section, we consider an image deblurring example where the desired image contains an inclusion whose degree of edge regularity can be characterized using a regularity parameter.  We consider training a DNN to learn \textit{both} the optimal TV regularization parameter and the regularity parameter of the inclusion.  It is natural to require the network to consider the dependency and coupling of these parameters.

\begin{figure}[t]
\centering
\begin{tabular}{cccc}
    \includegraphics[width=0.2\textwidth]{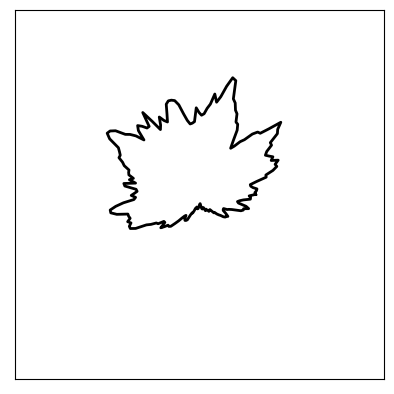} &
    \includegraphics[width=0.2\textwidth]{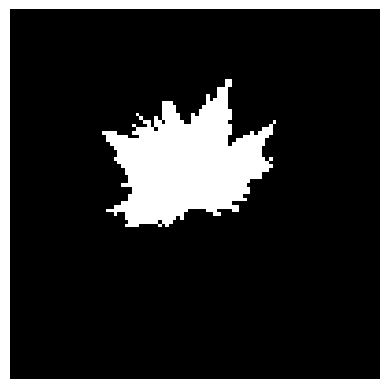} &
    \includegraphics[width=0.2\textwidth]{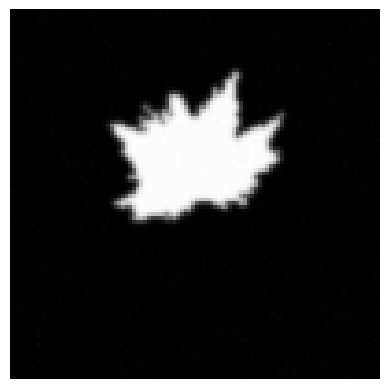} &
    \includegraphics[width=0.2\textwidth]{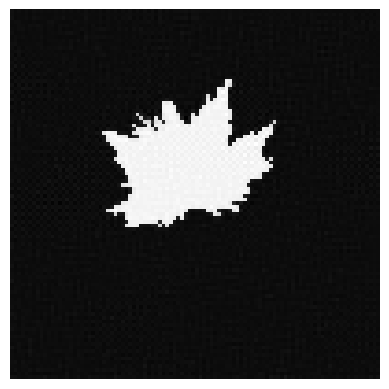} \\
    \includegraphics[width=0.2\textwidth]{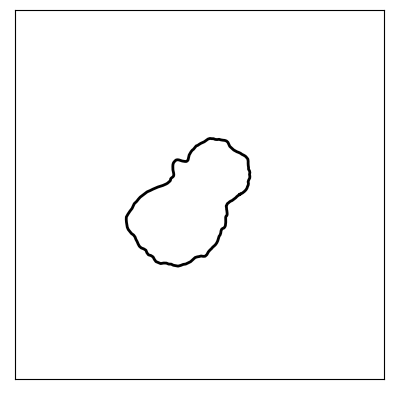} &
    \includegraphics[width=0.2\textwidth]{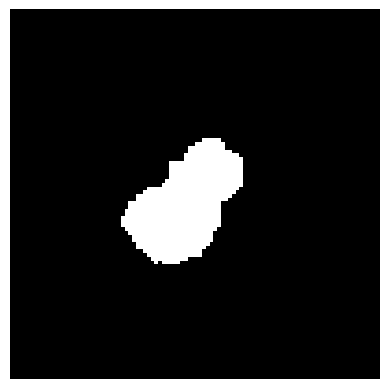} &
    \includegraphics[width=0.2\textwidth]{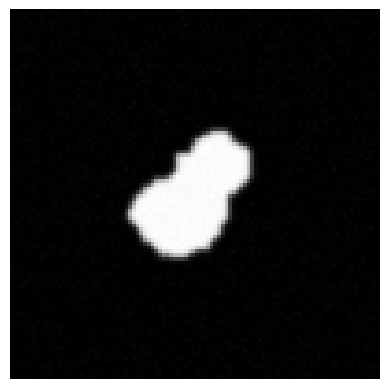} &
    \includegraphics[width=0.2\textwidth]{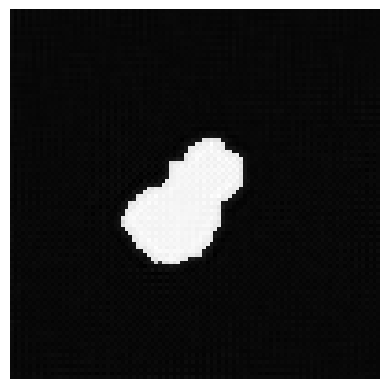} \\
    \includegraphics[width=0.2\textwidth]{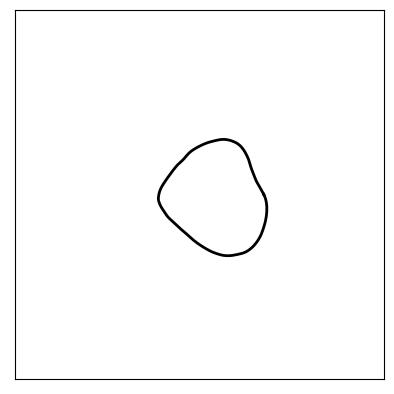} &
    \includegraphics[width=0.2\textwidth]{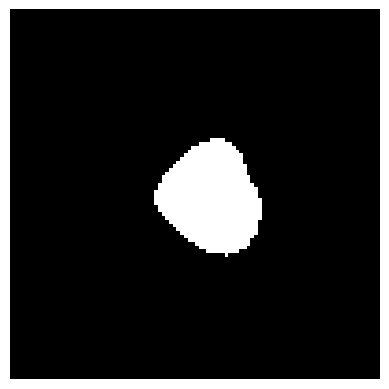} &
    \includegraphics[width=0.2\textwidth]{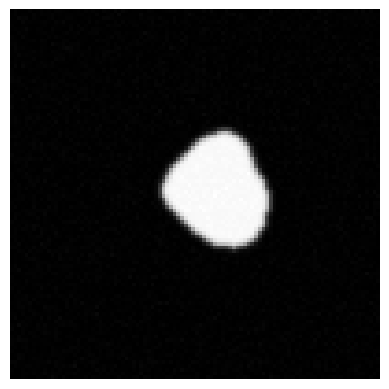} &
    \includegraphics[width=0.2\textwidth]{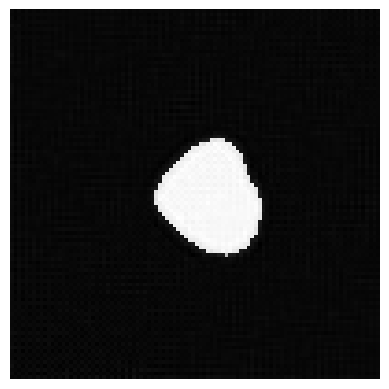} \\
    \makecell{ (a) outline of \\ the inclusion $r$ } & \makecell{ (b) ground \\ truth image $\mathbf x_{\text{true}}$} & \makecell{ (c) degraded \\ image $\mathbf b$ } & \makecell{(d) reconstructed \\ image $\widehat{\bfx}(\lambda_{\rm opt})$ }
\end{tabular}
\caption{Example images of star-shaped inclusions. In each row, we provide the outline of the inclusion, the true image, the blurred, noisy image, and the reconstruction using TV with the optimal regularization parameter. Regularity of $r$ is set to $\gamma = 1.2, 1.75$ and $2.5$, with noise level 1.79\%, 3.67\% and 4.07\%, for the first, second, and third row, respectively.}
\label{fig:blur_sample}
\end{figure}

We consider a discrete image deblurring example, where the vectorized blurred image contained in $\bfb$ can be modeled as~\eqref{eq:invprob} with $\bfA(\bfx_{\rm true}) = \bfA\bfx_{\rm true}$, see~\cite{hansen2006deblurring}.
Here, $\bfx_{\rm true}$ represents the vectorized desired image, $\bf A$ is a discretized linear blurring operator, and $\bfvarepsilon$ represents noise in the observed data. The simulated true images contain piecewise constant ``\emph{star-shaped}'' inclusions \cite{1930-8337_2016_4_1007,BuiThanh2014}, see Figure~\ref{fig:blur_sample} first column. Such inclusions are characterized by their center $\bfc_0$ and a radial function $r$. More precisely, let $\mathcal D \subset \bbR^2$ be the unit disk and $\tau:\bbR^2 \to \bbR$ be the continuous mapping from the Cartesian to angular polar coordinates. We define the region of inclusion to be
\begin{equation}
    \mathcal A(r,\bfc_0) = \left\{ \bfy\in \mathcal D : \| \bfy - \bfc_0 \|_2 \leq r(\tau( \bfy-\bfc_0)) \right\},
\end{equation}
where we set $\bfc_0$ to be the origin and radial function $r$ represents a 1-dimensional periodic log-Gaussian random field \cite{ibragimov2012gaussian} defined as
\begin{equation} \label{KL}
    r(\xi) = r_0 + c \exp\left( \frac{1}{\sqrt{\pi}} \sum_{i=1}^{\infty} \left(\frac{1}{i} \right)^\gamma \left(  X^1_i \cos(i\xi) + X^2_i \sin(i\xi) \right) \right).
\end{equation}
Here, we assume $X^1_i,X^2_i$ to be random normal $X^1_i,X^2_i \sim \mathcal N(0,1)$, $\gamma>1$ controls the regularity of the radial function, $r_0 > 0$ is the deterministic lower bound of the inclusion radius and $c>0$ is the amplification factor. We construct an infinite-dimensional image as
\begin{equation}
    \chi(r) = a^+ \mathbb 1_{ \mathcal A } + a^- \mathbb 1_{ \mathcal D \backslash \mathcal A },
\end{equation}
where $\mathbb 1$ is the indicator function, $a^+ = 1$ and $a^-=0$. Discretizing $\chi$ into pixels and reshaping it into a vector yields the ground truth discrete image $\bfx_{\rm true}$. For this test case, we fix the size of the images to $100\times 100$ pixels. We generate a dataset of $J = 20,\!000$ true images $x_{\text{true}}^{(j)}$ by selecting the inclusion regularity parameter $\gamma$ uniformly from the interval $\gamma \in [1.25, 2.5]$, setting the application factor $c=0.25\exp(0.2)$ and setting the deterministic minimum value to $r_0 = 0.2/\sqrt{2\pi}$. We truncate the sum in the log-Gaussian random field after 100 terms.
Furthermore, we simulate the corresponding observed images $\mathbf b^{(j)} = \mathbf A \mathbf x_{\rm true}^{(j)} + \bfvarepsilon^{(j)}$, where $\mathbf A$ comprises a 2-dimensional discrete Gaussian blur with a standard deviation of $\sigma_{\kappa} = 1$ and a stencil of the size $5\times 5$ pixels and $\bfvarepsilon^{(j)}$ is white noise where the noise level is selected uniformly between $0.1\%$ and $5\%$.  Three sample inclusions corresponding to different choices of $\gamma$ are provided in Figure~\ref{fig:blur_sample}, along with the true and observed images.

Regularization is an essential step in solving the image deblurring
problem, and many choices for the regularization term $\cal R(\bf x)$ have been considered \cite{hansen2006deblurring,5701777,Huang2008,7448477,6909767}. We consider TV regularization~\eqref{eq:TV} as discussed earlier, since it provides an excellent choice for deblurring images with piecewise smooth components.
Larger values of the partial derivatives are only allowed in certain regions in the image \cite{hansen2006deblurring}, e.g., near edges and discontinuities.
Thus, an optimal choice of the regularization parameter $\lambda$ in \eqref{eq:TV} depends on the edge properties of the particular image. For example, notice that for smaller $\gamma$ the true star-shaped inclusion in $\mathbf x_{\rm true}$ contains a longer boundary (in fact, the length of the boundary tends to $\infty$ as $\gamma \to 1$.) In this case, the TV regularization term~\eqref{eq:TV} will have a larger contribution to the minimization problem~\eqref{eq:loss}.
The dependency of $\lambda_{\rm opt}$ on $\gamma$ is in general nonlinear.  We remark that an accurate prediction of $\gamma$ can have significance beyond its impact on the optimal regularization parameter.
For example, in some applications the prediction of the outline of inclusions in a degraded image can be used to differentiate cancerous versus non-cancerous tissues. Furthermore, predicting the right regularity of the outline, i.e., $\gamma$, can significantly enhance the uncertainty quantification of such a prediction.

Next we describe the learning process. For each blurred image, we solve~\eqref{eq:lambdatrue} to obtain $\lambda_{\rm opt}^{(j)}$, using the split-Bregman approach \cite{rudin1992nonlinear} to find $\widehat{\bfx}(\lambda)$ that solves \eqref{eq:1lambda} with~\eqref{eq:TV}. We denote the corresponding reconstruction by $\bfx_{\rm opt}^{(j)} = \widehat\bfx(\lambda_{\rm opt}^{(j)})$, see Figure~\ref{fig:blur_sample}.
Thus, the dataset consists of $\left\{ \mathbf{b}^{(j)}, \mathbf{x}_{\rm true}^{(j)}, \lambda_{\rm opt}^{(j)}, \gamma_{\rm true}^{(j)} \right\}$. We construct a DNN to predict the regularity parameter $\gamma_{\rm true}$ and the optimal regularization parameter $\lambda_{\rm opt}$, simultaneously. The DNN comprises of 3 convolutional sub-networks and 2 types of output networks for $\gamma_{\rm true}$ and $\lambda_{\rm opt}$, respectively. Each convolutional sub-network consists of a 2-dimensional convolution layer, a 2-dimensional batch normalization, a ReLU activation function, and a 2-dimensional max-pool layer. We consider an output layer with a single ReLU layer for $\gamma_{\rm true}$ (out\_1 in Table~\ref{tab:CNN_star}) and an output layer for $\lambda_{\text{opt}}$ with 3 hidden layers (out\_2 in Table \ref{tab:CNN_star}). Since the dependency of $\lambda_{\text{opt}}$ on $\gamma$ is nonlinear in general, the output network for $\lambda_{\text{opt}}$ is chosen deeper than the output network for $\gamma$ to capture the extra complexity. We summarize the architecture of this DNN in Table~\ref{tab:CNN_star}.

\begin{table}[tbhp]
{\footnotesize
  \caption{Description of DNN for image deblurring with star-shaped inclusions example.}\label{tab:CNN_star}
\begin{center}
\begin{tabular}{|c|c|c|c| } \hline
Layer name & number of layers & output size & layer type \\ \hline
conv$_1$ & - & 50$\times$ 50 with 16 channels & convolutional with kernel size:$5\times 5$ \\ \hline
conv$_2$ & - & 25$\times$ 25 with 32 channels & convolutional with kernel size:$5\times 5$ \\ \hline
conv$_3$ & - & 12$\times$ 12 with 64 channels & convolutional with kernel size:$5\times 5$ \\ \hline
out$_1$ & 1 & 1 & feed-forward \\ \hline
out$_2$ & 3 & 1 & feed-forward \\ \hline
\end{tabular}
\end{center}
}
\end{table}

Let $\widehat\bfPhi_{\rm conv}$ denote the convolutional part of the network and $\widehat\bfPhi_{\gamma}$ and $\widehat\bfPhi_{\lambda}$ denote the output networks corresponding to $\gamma_{\rm true}$ and $\lambda_{\rm opt}$, respectively. Since the DNN produces multiple outputs, we train it in 2 separate stages. In the first stage we train the network $\widehat\bfPhi_{\gamma} \circ \widehat\bfPhi_{\rm conv}$ which maps $\mathbf b$ onto $\gamma$. In the second stage we fix the network parameters in $\widehat\bfPhi_{\rm conv}$, denote the fixed network by $\widehat{ \bfPhi}_{\rm conv}^{\gamma}$, and train the network $\widehat{ \bfPhi}_{\gamma} \circ \widehat{ \bfPhi }_{\rm conv}^{\gamma}$. The two stages for updating the network can be carried out for a single batch of data, or after a complete training of each network. The cost function in the first stage is chosen to be~\eqref{eq:dnnlearn} where $\lambda_{\rm opt}^{(j)}$ is replace by $\gamma_{\rm true}^{(j)}$.

Recall that the convolutional part of the network extracts information in an image that is required in generating the output. The 2-stage method in training the network assumes that the necessary information required for reconstructing $\gamma_{\rm true}$ is the same as that needed for the reconstruction of $\lambda_{\rm opt}$. Furthermore, this technique conserves the order of dependency between the parameters, i.e., from $\gamma_{\rm true}$ to $\lambda_{\rm opt}$.

\begin{figure}
\centering
    \includegraphics[width=0.45\textwidth]{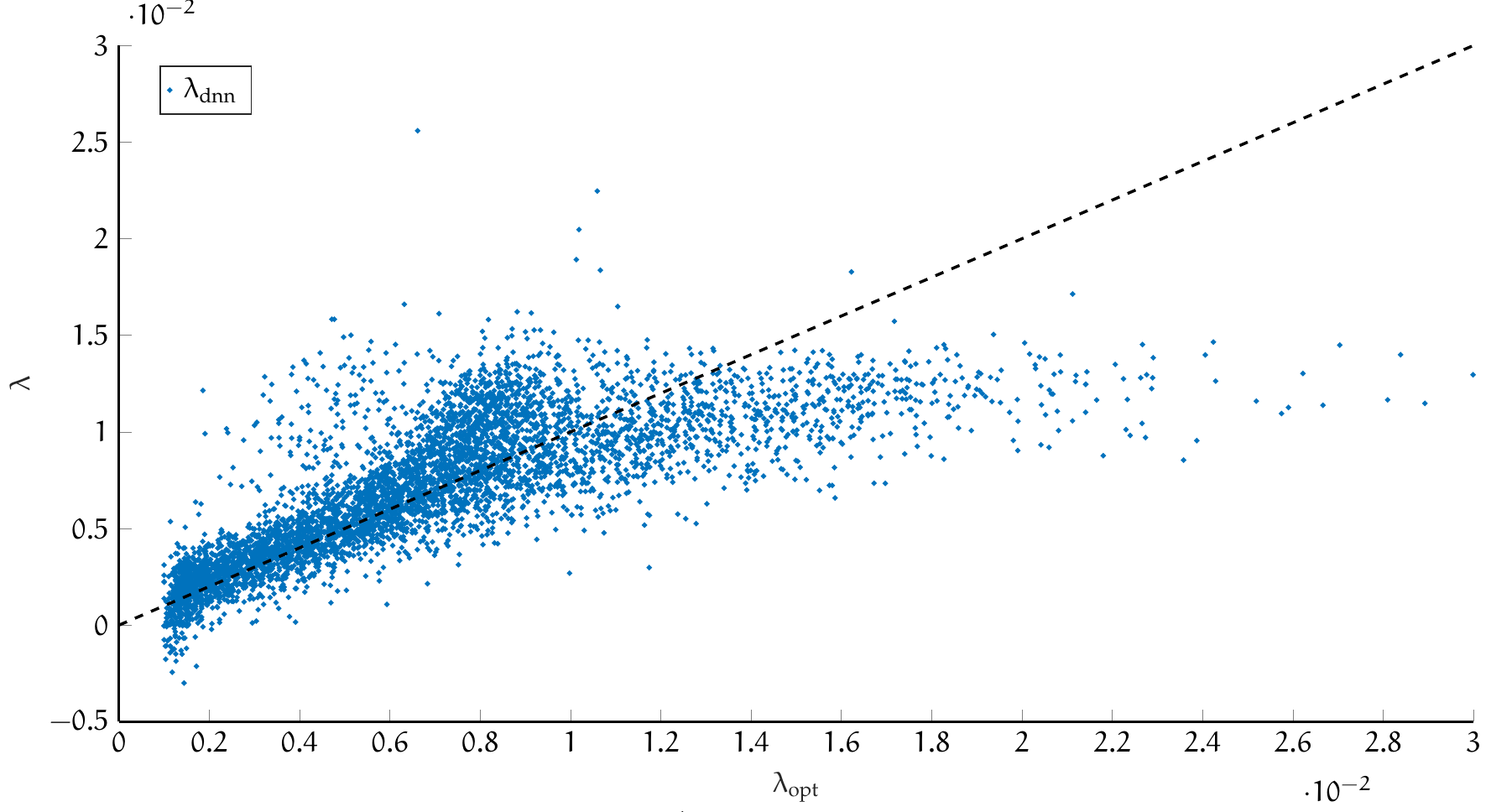}
    \includegraphics[width=0.45\textwidth]{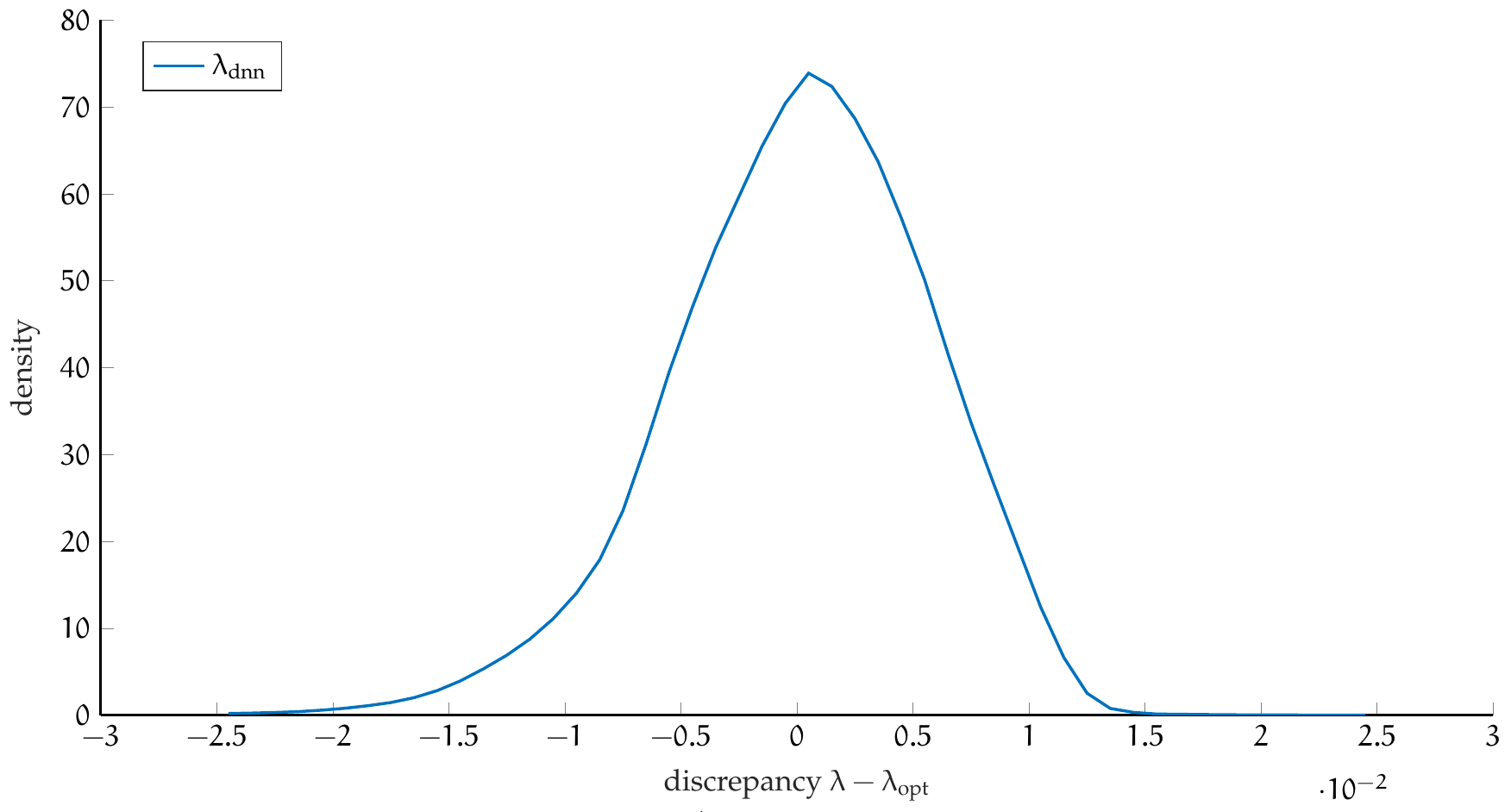}
\caption{Prediction of regularization parameter $\lambda_{\rm opt}$ for the image deblurring with inclusions example with uniformly distributed $\gamma_{\rm true}\in [1.2,2.5]$. (left) Scatter plot of network predicted regularization parameter $\lambda_{\rm dnn}$ versus the optimal regularization parameter $\lambda_{\rm opt}$.
(right) Probability density for the discrepancies between the network predicted regularization parameter and the optimal regularization parameter.}
\label{fig:star_net}
\end{figure}

\begin{figure}
\centering
    \includegraphics[width=0.7\textwidth]{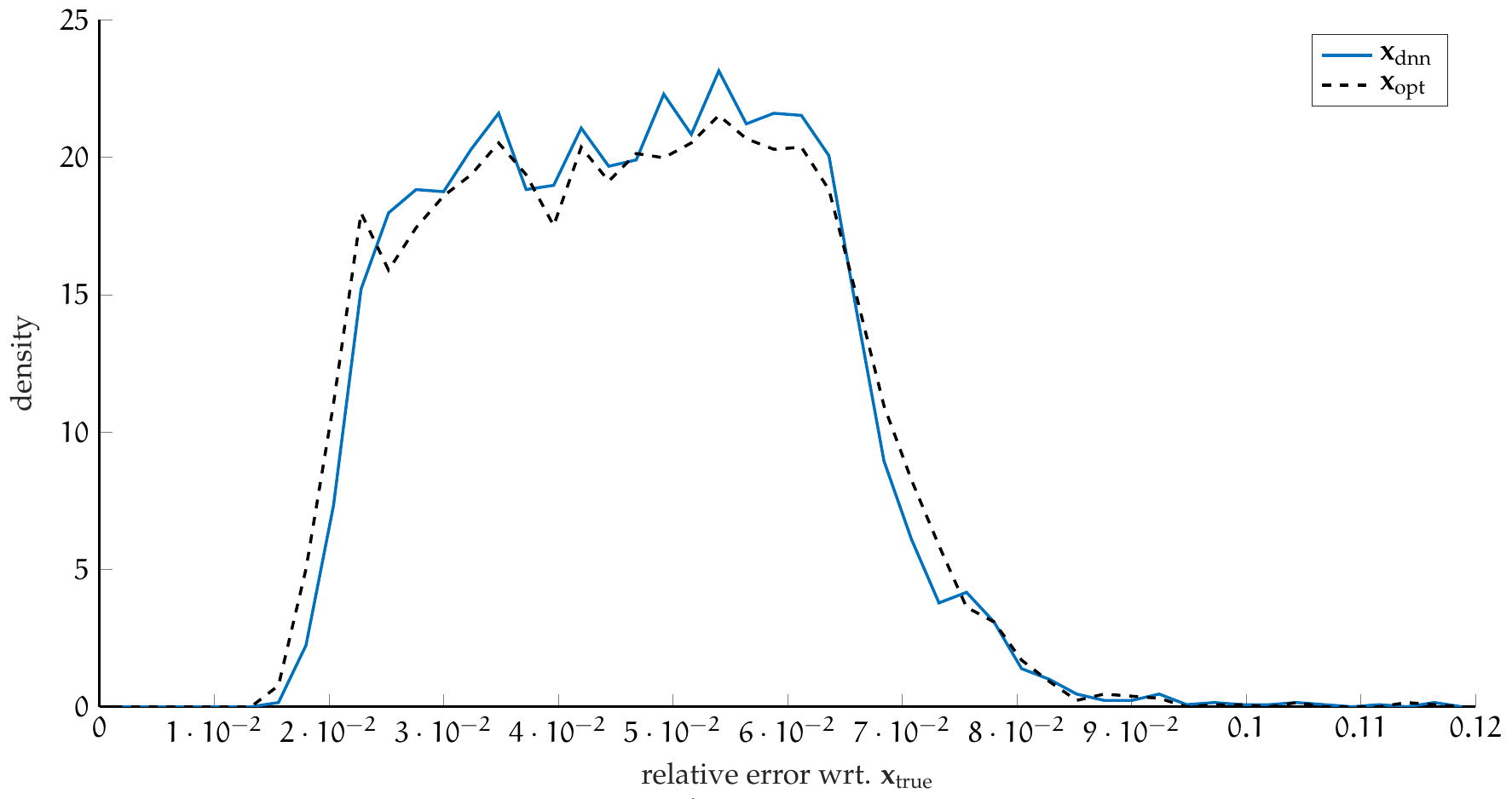}
\caption{Probability distribution of relative reconstruction error norms computed with respect to the $\ell^1$-norm for the image deblurring with inclusions example.
}
\label{fig:star_hist}
\end{figure}

To train the network, we split the dataset into a training set of $15,000$ data points, and a validation set of $5,000$ data points. We utilize the ADAM optimizer with a dynamic learning rate chosen in the interval $[10^{-5}, 10^{-1}]$ with a batch size of $2^{10}$ data points. The network is trained for $10^4$ epochs. The performance of the method is evaluated on the validation set.  In Figure~\ref{fig:star_net}, we summarize the performance of the network in predicting the optimal regularization parameter. The scatter plot in the left panel indicates a strong correlation between the optimal and the DNN predicted regularization parameter, and the plot in the right panel indicates relatively small discrepancies.
The performance of the predicted regularization parameter in terms of the reconstructed image can be found in Figure \ref{fig:star_hist}, where we provide the probability densities for the relative reconstruction error norms compared to $\bfx_{\rm true}$ and calculated with respect to the $\ell^1$-norm.
We observe an excellent match between the distribution of the error norms for the image reconstructed by $\lambda_{\rm opt}$ and for the image reconstruction by $\lambda_{\rm dnn}$. The authors found comparable results when an $\ell^2$-norm is utilized for the relative reconstruction errors.

\begin{figure}
\centering
    \includegraphics[width=0.45\textwidth]{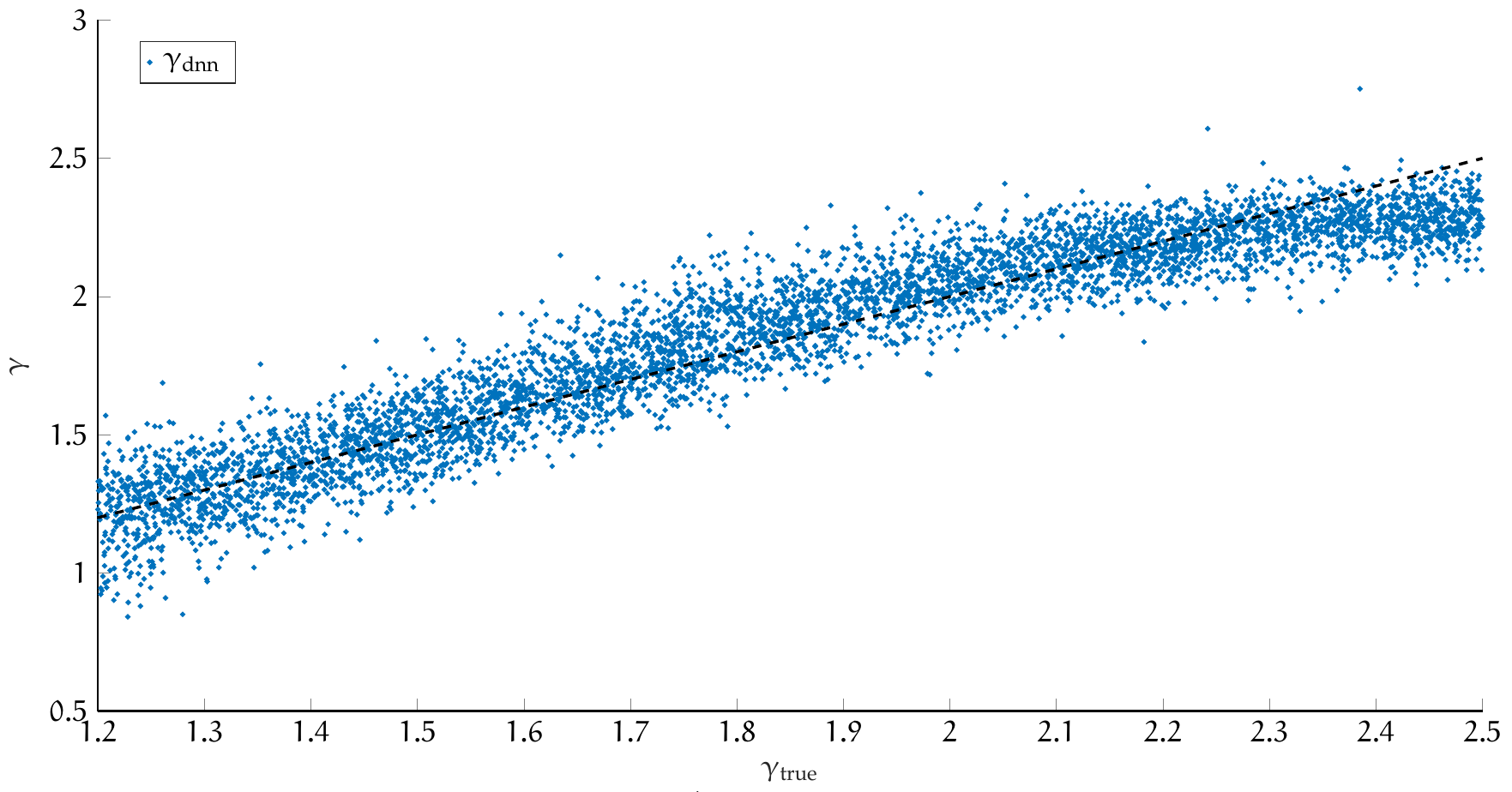}
    \includegraphics[width=0.45\textwidth]{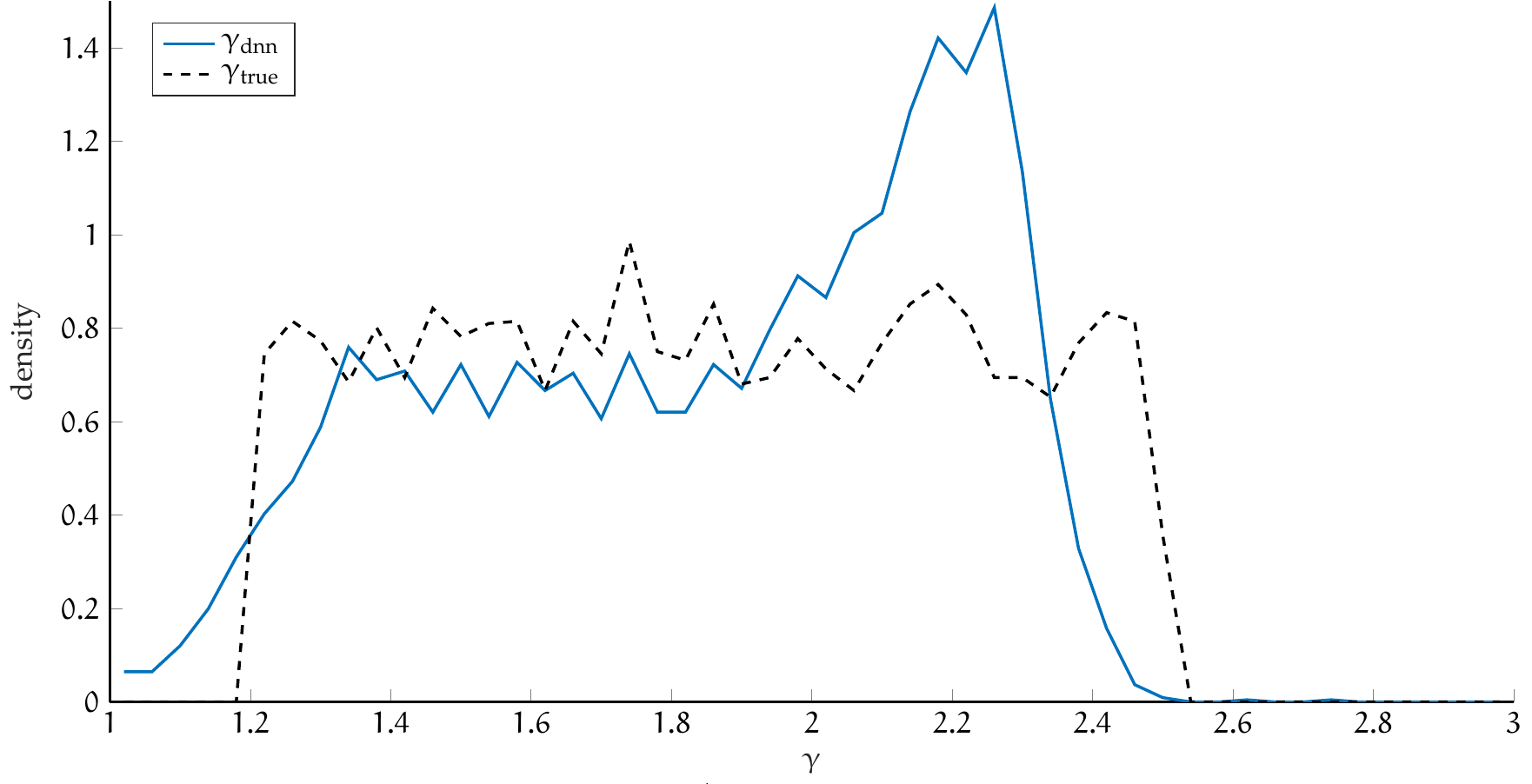}
\caption{Prediction of regularity parameter $\gamma_{\rm true}$ for the image deblurring with inclusions example with uniformly distributed $\gamma_{\rm true}\in [1.2,2.5]$. (left) Scatter plot of network predicted regularity parameter $\gamma_{\rm dnn}$ versus the true regularity parameter $\gamma_{\rm true}$.
(right) Probability densities of $\gamma_{\rm dnn}$  and $\gamma_{\rm true}$.}
\label{fig:star_gamma}
\end{figure}

Finally, we investigate the performance of the DNN in predicting $\gamma$, the parameter defining the regularity of the star-shaped inclusion. The scatter plot in the left panel of Figure~\ref{fig:star_gamma} shows high correlation between the true regularity parameter $\gamma_{\rm true}$ and the DNN predicted parameter $\gamma_{\rm dnn}$. The plot in the right panel of Figure \ref{fig:star_gamma} provides probability densities of $\gamma_{\rm true}$ and $\gamma_{\rm dnn}$.
We notice larger errors in the prediction for larger values of $\gamma$. This can be due to the low resolution of images where the smoothness information is lost in the discretization. Recall that $\widehat\bfPhi_{\rm conv}$ only extracts information in an image that is needed to predict $\gamma$. Figure~\ref{fig:star_hist} validates that this information is sufficient for the prediction of $\lambda_{\rm opt}$.

\subsection{Learning the stopping iteration for iterative regularization}
\label{sec:iter}
In this example, we train a convolutional neural network to learn the mapping from observation to optimal stopping iteration. We consider a linear inverse diffusion example described in \cite{gazzola2019ir,min2013inverse} where the goal is to determine an initial function, given measurements obtained at some later time.  The solution is represented on a finite-element mesh and the forward computation involves the solution of a time-dependent PDE. The underlying problem is a 2D diffusion problem in the domain $[0, T]\times[0, 1]\times[0, 1]$ in which the solution
$x$ satisfies
\begin{equation}
\label{eq:pde}
    \frac{\partial x}{\partial t} = \nabla^2 x
\end{equation}
with homogeneous Neumann boundary conditions and a smooth function $x_0$ as initial condition at time $t = 0$. The forward problem maps $x_0$ to the solution $x_T$ at time $t = T$, and the inverse problem is then to reconstruct the initial condition from observations of $x_T$.
We discretize the function $x$ on a uniform finite-element mesh with $2(\sqrt{n}-1)^2$ triangular elements, where the domain is an $(\sqrt{n}-1) \times (\sqrt{n} - 1)$ pixel grid with two triangular elements in each pixel. Then, vector $\bfx \in \bbR^n$ contains the $n$ values at the corners of the elements. The forward computation is the numerical solution of the PDE \eqref{eq:pde} using the Crank-Nicolson-Galerkin finite-element method, and the discretized forward process is represented $\bfA \in \bbR^{n\times n}$, see \cite{gazzola2019ir}.

\begin{figure}[bthp]
\begin{center}
 \includegraphics[width=0.35\textwidth]{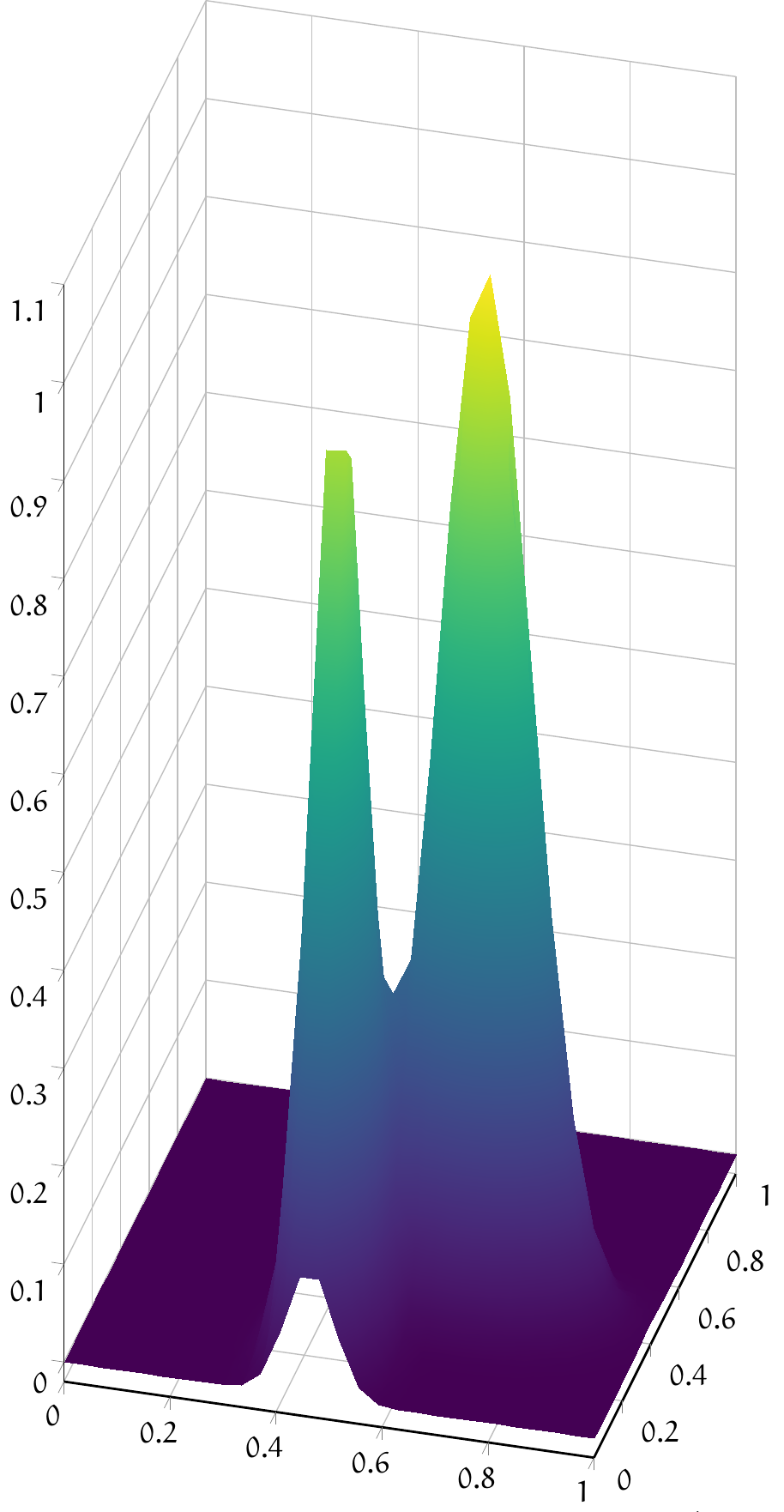}
 \includegraphics[width=0.35\textwidth]{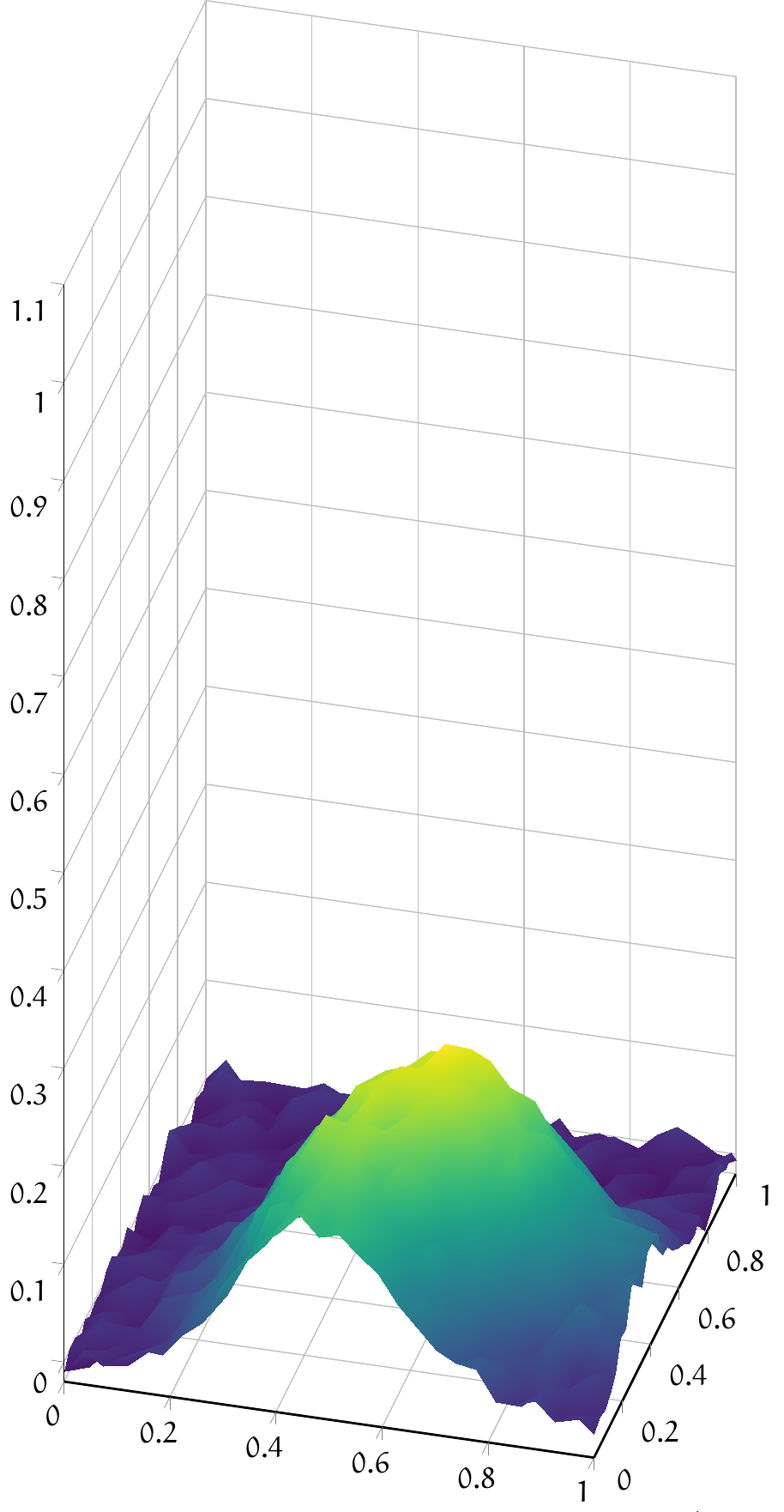}
 \end{center}
\caption{2D inverse diffusion problem: For $n = 28$, we provide one example of the true solution $\bfx_{\rm true}$ corresponding to the initial function $x_0$ on the left and the corresponding observed data $\bfb$ at time $T = 0.01$ on the right.}
\label{fig:recon1}
\end{figure}

We generate a data set containing initializations $\bfx_{\rm true}^{(j)} = \bfx_0^{(j)}$ for $j=1, \ldots, 12,\!500$ where $n = 784 = 28 \times 28$.  Each initialization is generated as
\begin{equation}
    x_0(\bfxi) = a \psi(\bfxi,\bfc_1, \bfnu_1) + \psi(\bfxi,\bfc_2, \bfnu_2),
\end{equation}
where $\bfxi\in \bbR^2$ represents the spatial location, $\psi(\bfxi, \bfc, \bfnu) = \e^{-(\bfxi- \bfc)\t \diag{\bfnu} (\bfxi- \bfc)}$
where the amplitude $a=0.7|\zeta|$ where $\zeta \sim \calN(0,1)$, the components of the centers $\bfc_1, \bfc_2\in \bbR^2$ are uniformly randomly selected from the interval $[0.1, 0.9]$, and the components of vectors $\bfnu_1, \bfnu_2\in \bbR^2$ are uniformly randomly selected from the interval $[5\cdot10^{-2}, 2\cdot 10^{-1}]$.
The $j$-th observation is generated as
\begin{equation}
\label{eq:diffusionmodel}
    \bfb^{(j)} = \bfA \bfx_{\rm true}^{(j)} + \bfvarepsilon^{(j)}
\end{equation}
where the noise level is uniformly randomly selected from the interval $[10^{-5}, 5\cdot10^{-1}]$. An example of one initialization along with the corresponding noisy observed data $\bfb$ is provided in Figure~\ref{fig:recon1}.  The noise level for this example is $0.0974$.

Next, we consider the reconstruction process for the 2D inverse diffusion problem. Although many iterative projection methods could be used here, we consider the range-restricted GMRES (RRGMRES) method, which does not require the transpose operation. For each sample, we run RRGMRES to compute the corresponding optimal stopping iteration, i.e., the stopping iteration that corresponds to the smallest relative reconstruction error,
\begin{equation}
    k_{\rm opt}^{(j)} = \argmin_{k\in \bbN}  \ \|\bfx_k - \bfx_{\rm true} \|
\end{equation}
where $\bfx_k$ is the $k$-th iterate of the RRGMRES approach applied to \eqref{eq:diffusionmodel}.

\begin{figure}[bthp]
\begin{center}
\includegraphics[width=0.7\textwidth]{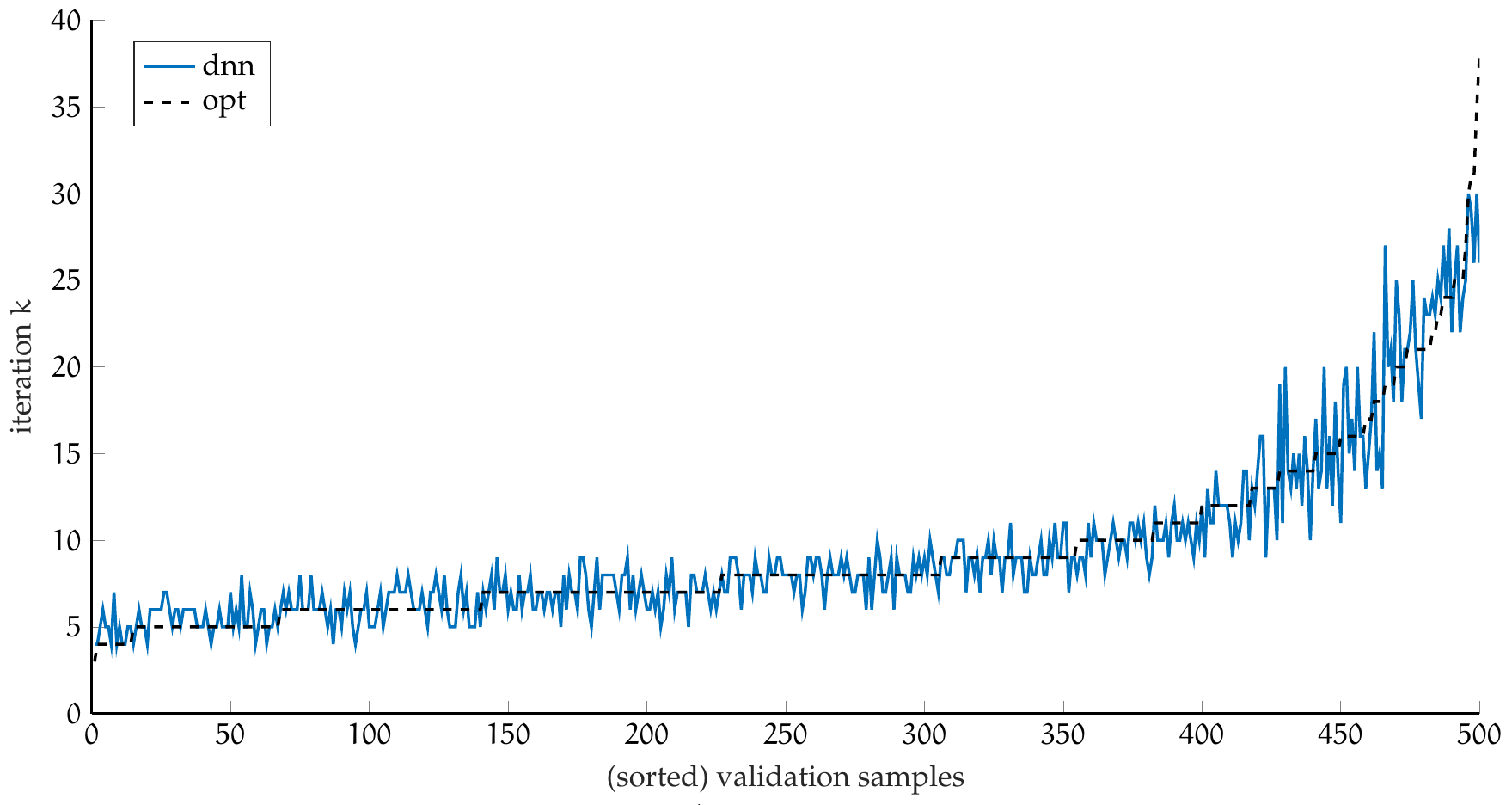}
\end{center}
\caption{Optimal and DNN predicted stopping iteration for each sample of the validation set.}
\label{fig:diff_iterations_all}
\end{figure}

Now we have a data set containing $12,\!500$ pairs, $\left\{\bfb^{(j)}, k_{\rm opt}^{(j)}\right\}$.  We split the data into $12,\!000$ samples for the training data and $500$ samples in the validation data. Using the training data, we employ a convolutional neural network with four convolutional layers. Each convolutional layer consists of $3\times 3$ filters with the following number of channels $8, 16, 32$ and $32$ respectively and a bias term for each. The convolutional layers are padded, and the stride is set to~1.
To reduce the dimensionality of the neural network we use an average pooling of $2\times 2$. We establish one $20$\% dropout layer followed by one dimensional output layer $1\times 1,\!568$, plus bias term.
Each hidden layer has a ReLU activation function. To estimate $\bftheta$ we utilize the stochastic gradient descent with momentum method with a learning rate of $10^{-4}$, while the batch size is set to $128$. We learn for $50$ epochs.
Although the stopping iteration must be a whole number that is greater than $1$, we used a regression loss output layer and just rounded the outputs of the DNN for prediction. The regression loss assumes that the distribution of errors is Gaussian which is not the case with an integer output. We remark that a more suitable loss function (e.g., a Poisson loss or negative binomial loss) could be used.
For the $j$-th sample, the DNN predicted stopping iteration is denoted by $k_{\rm dnn}^{(j)}.$
In Figure~\ref{fig:diff_iterations_all}, we plot the optimal iteration $k_{\rm opt}$ along with the DNN predicted stopping iteration $k_{\rm dnn}$ per  (sorted) validation sample. Notice that the predicted stopping iteration via the learned DNN network is close to the optimal stopping iteration.

For the validation set, we provide comparisons to results using the DP to estimate the stopping iteration.  For these results, for each sample we estimate the noise level $\eta$ from the data $\bfb^{(j)}$ using a wavelet noise estimator \cite{donoho1995noising} and select the iteration $k_{\rm dp}^{(j)}$ such that $\norm[2]{\bfb^{(j)} - \bfA \bfx_k^{(j)} }/\norm[2]{\bfb^{(j)}} \leq \eta$
where $\bfx_k^{(j)}$ is the $k$-th iterate of RRGMRES and the safety factor $\eta = 1.01$ was suggested in \cite{gazzola2019ir}. We found that there were some examples in the validation set where the DP failed, resulting in very large reconstruction errors.  Of the $500$ validation examples, there were $25$ examples where DP used to compute a stopping iteration resulted in relative reconstruction error norms above $2$. For visualization purposes, these are not provided in the following results.

\begin{figure}[bthp]
    \centering
    \includegraphics[width=0.45\textwidth]{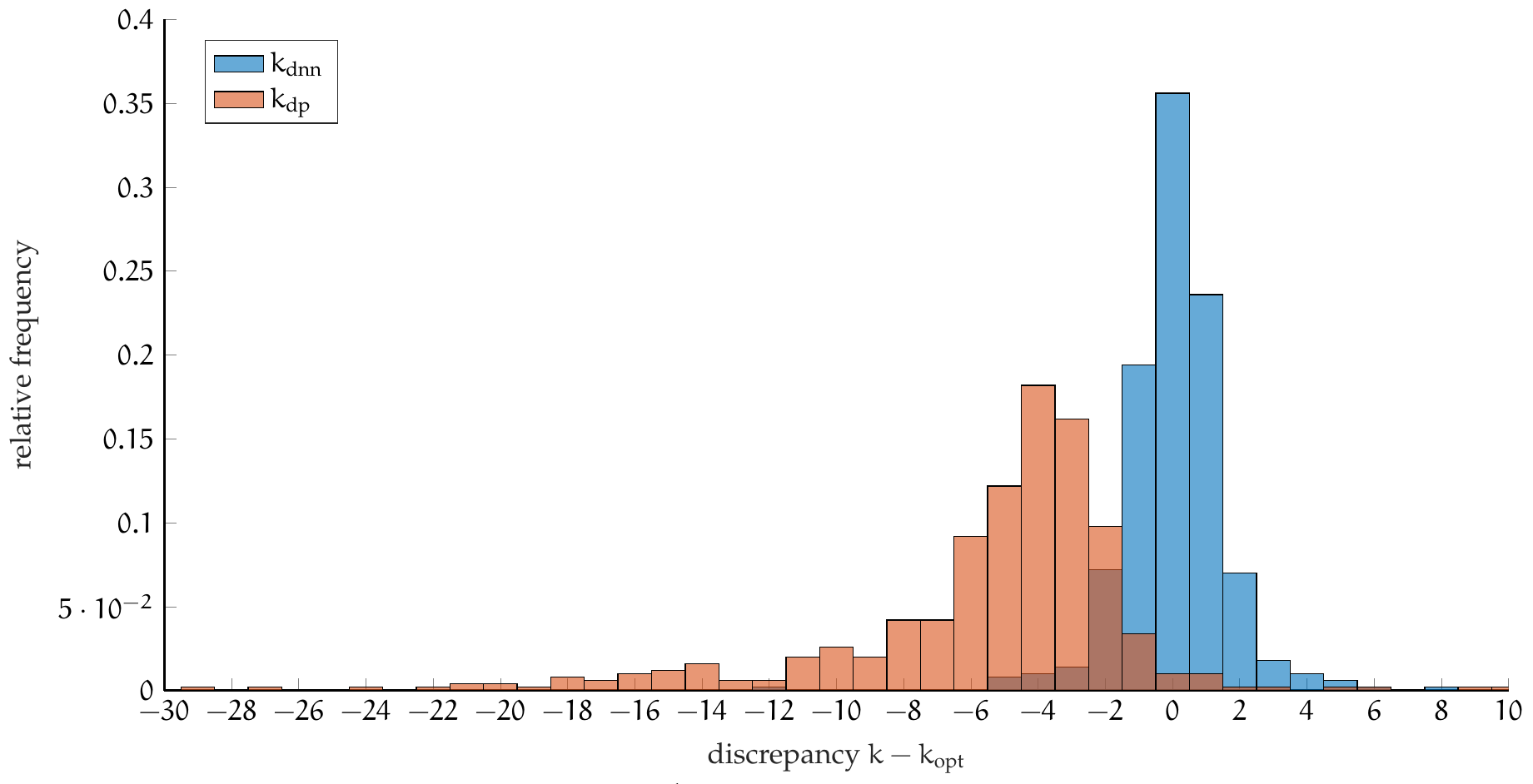}
    \includegraphics[width=0.45\textwidth]{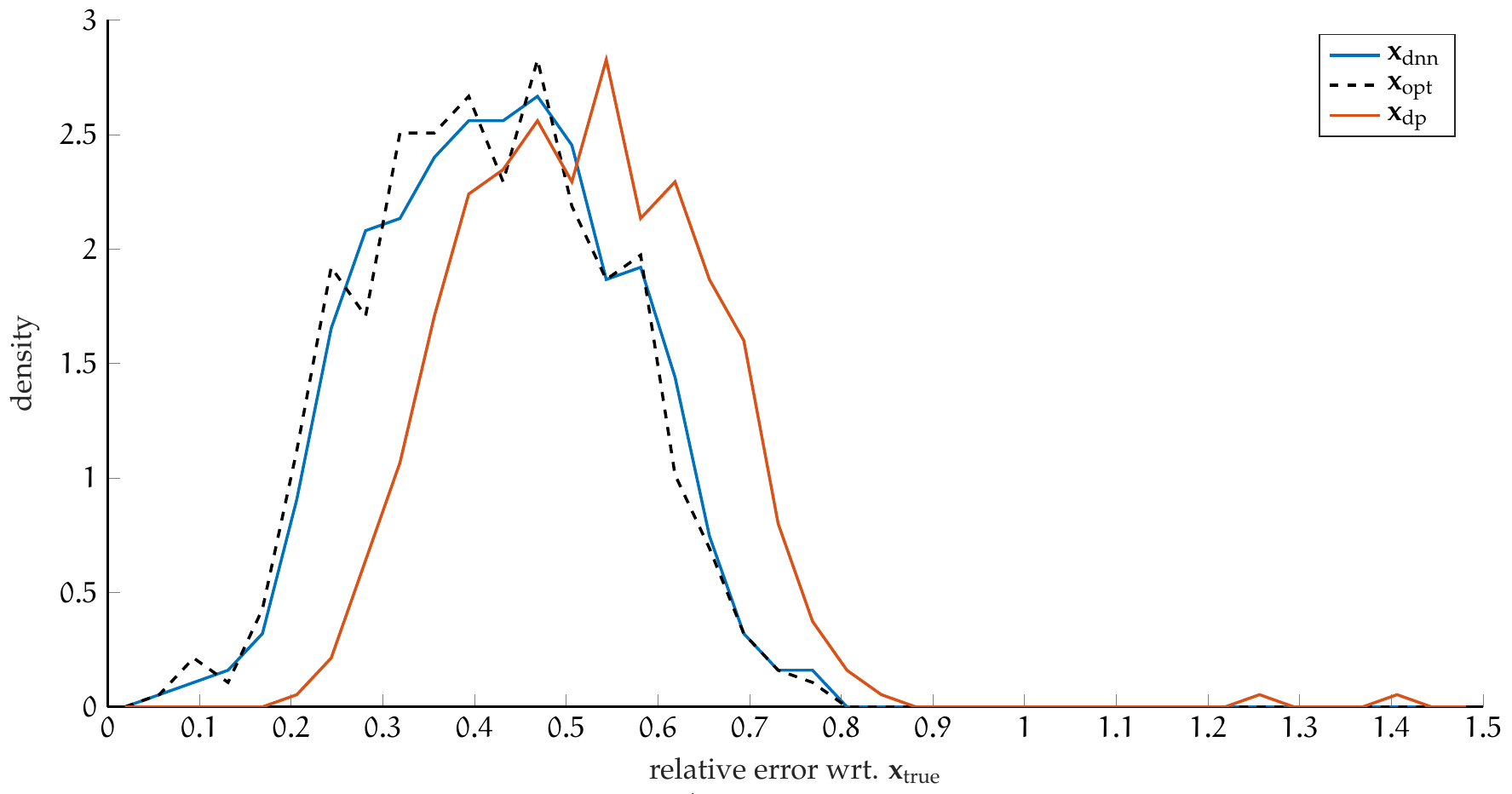}
    \caption{The left panel depicts the distribution of the discrepancy between the estimated stopping iteration $k$ for DNN and DP and the optimal stopping iteration $k_{\rm opt}$ for $500$ validation data. In the right panel, we provide the corresponding densities of relative reconstruction error norms.}\label{fig:diff_all}
\end{figure}

In left panel of Figure~\ref{fig:diff_all} we provide the distribution of the discrepancies between the DNN predicted stopping iteration and the optimal stopping iteration, $k_{\rm dnn}^{(j)} - k_{\rm opt}^{(j)}$. Notice that the distribution of discrepancies is centered around zero.  For comparison, we also provide the distribution of discrepancies for the DP, $k_{\rm dp}^{(j)} - k_{\rm opt}^{(j)}$.  We observe that the DP often underestimates the optimal stopping iteration. In the right panel of Figure~\ref{fig:diff_all} we provide the distribution of the relative reconstruction error norms with respect to $\bfx_{\rm true}$, i.e., $\norm[2]{\bfx_k - \bfx_{\rm true}}/\norm[2]{\bfx_{\rm true}}$ where $\bfx_k$ are reconstructions at stopping iterations $ k_{\rm dnn}^{(j)}$, $k_{\rm opt}^{(j)}$, and $k_{\rm dp}^{(j)}$. We observe that the DNN predicted stopping iterations result in relative reconstruction errors that are very close to those at the optimal stopping iteration.

\begin{figure}[bthp]
\begin{center}
\includegraphics[width=0.7\textwidth]{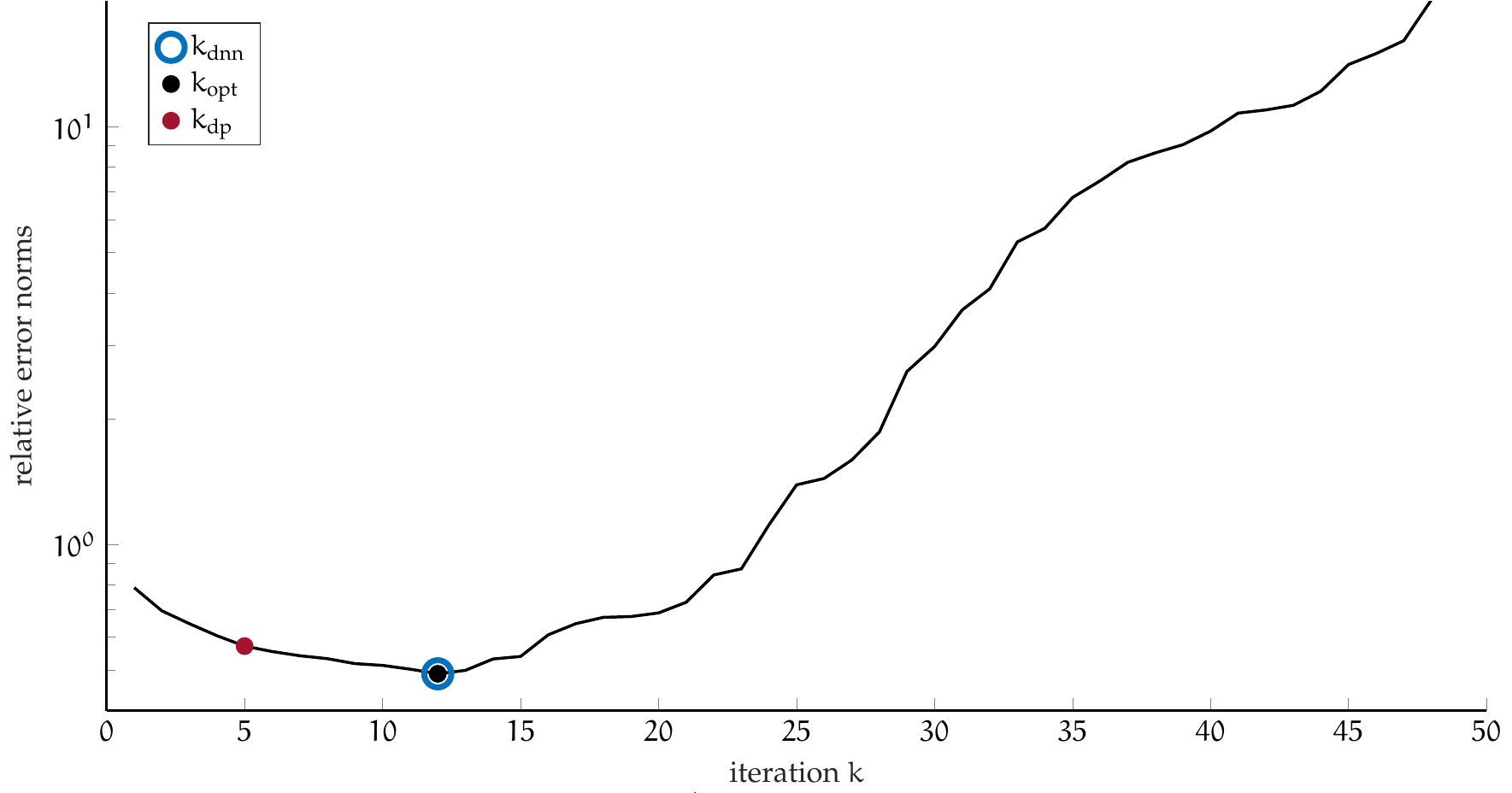}
\end{center}
\caption{For the example in Figure~\ref{fig:recon1}, we provide the relative reconstruction error norms per iteration of the RRGMRES method.  The markers correspond to the stopping iteration that is predicted via the learned DNN, the optimal stopping iteration, and the DP-selected stopping iteration.}
\label{fig:relerrors}
\end{figure}

\begin{figure}[bthp]
  \begin{center}
  \includegraphics[width=0.3\textwidth]{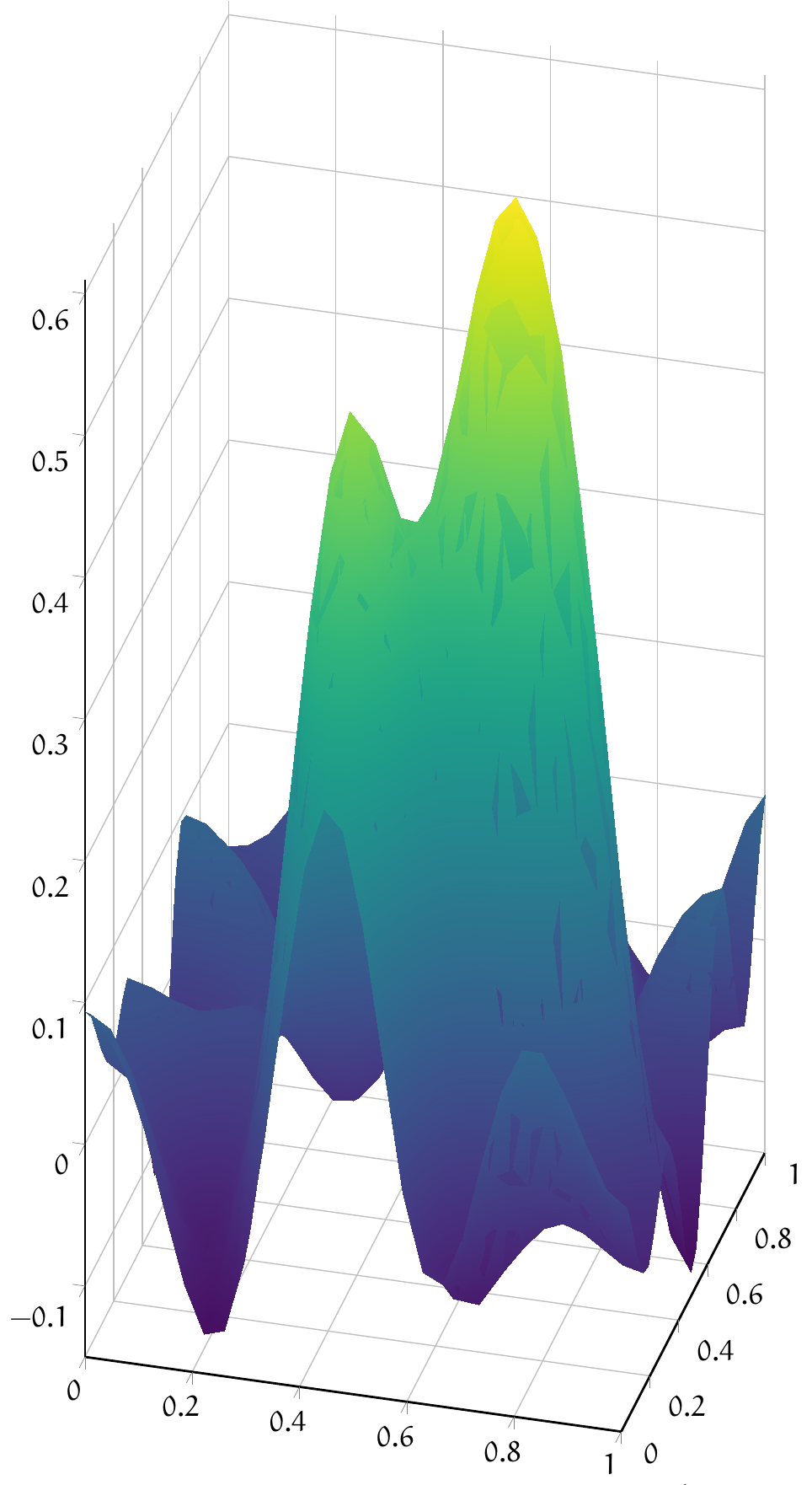}
  \includegraphics[width=0.3\textwidth]{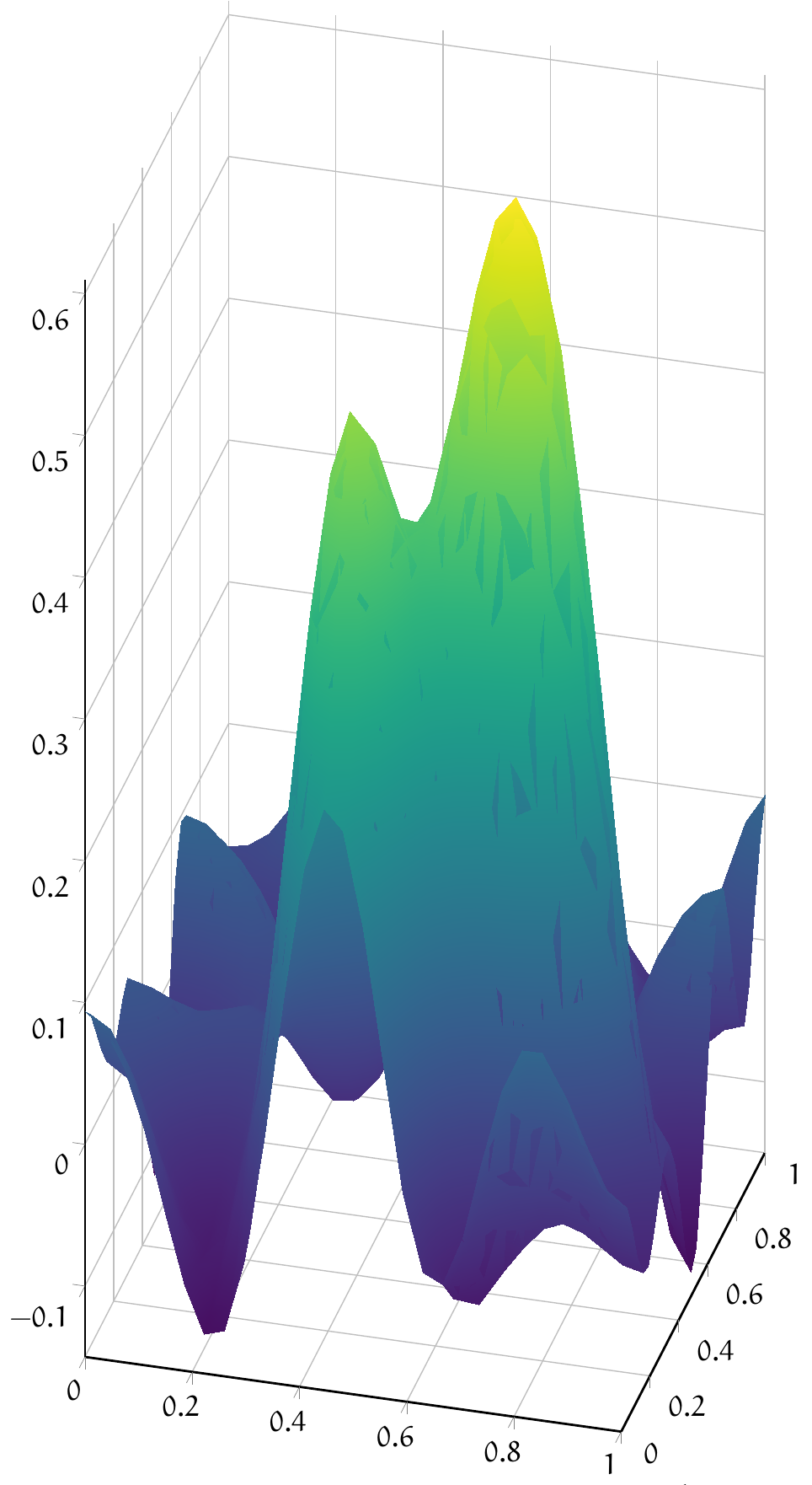}
  \includegraphics[width=0.3\textwidth]{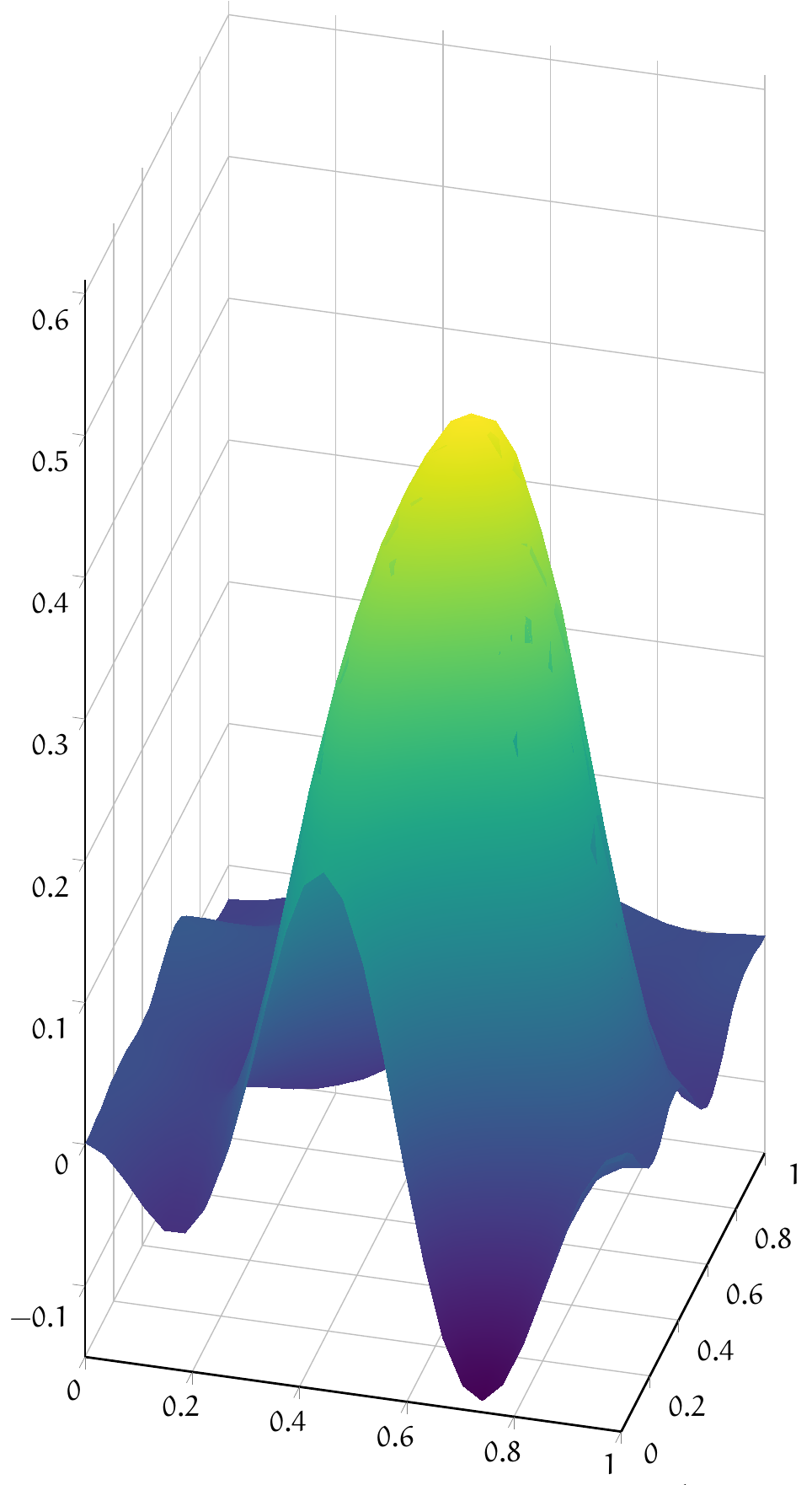}
  \end{center}
\caption{Reconstructions obtained using RRGMRES for the 2D inverse diffusion example in Figure~\ref{fig:recon1}, where the stopping iteration was determined via DNN, optimal, and DP.}
\label{fig:recon2}
\end{figure}

Last, we provide reconstructions corresponding to one sample from the validation set, where the true and observed signals are provided in Figure~\ref{fig:recon1}.  In Figure~\ref{fig:relerrors}, we provide the relative reconstruction errors per iteration of RRGMRES.  The optimal stopping iteration (corresponding to the minimizer of the relative reconstruction error norms is $12$ and is marked in black.  The DNN predicted stopping iteration was also $12$ and is marked by the blue circle.  The DP stopping iteration was $5$ and is denoted in red.  Corresponding reconstructions are provided in Figure~\ref{fig:recon2}, where it is evident that with DP, the reconstruction is too smooth and unable to resolve the two peaks in the initialization.

\section{Conclusions}
\label{sec:conclusions}
In this paper, we propose a new approach that uses DNNs for computing regularization parameters for inverse problems.  Using training data, we learn a neural network that can approximate the mapping from observation data to regularization parameters.  We consider various types of regularization including Tikhonov, total variation, and iterative regularization.  We also showed that this approach can be used to estimate multiple parameters (e.g., regularity of edges of inclusions and the regularization parameter).
We showed that DNN learned regularization parameters can be more accurate than traditional methods (not just in estimating the optimal regularization parameter but also in the corresponding reconstruction) and can be obtained much more efficiently in an online phase (requiring only a forward propagation through the network).  Although the proposed approach bears some similarity to existing OED approaches since the main computational costs are shifted to the offline phase, the DNN approach exhibits better generalizability since the computed regularization parameters are tailored to the specific data. Our results demonstrate that the mapping from the observation $\bfb$ to the regularization parameters $\bflambda$ can be well-approximated by a neural network.  We observed that our approach is flexible with regards to the specific design of the network and that despite the large dimension of the network input $\bfb$, not a significant amount of training data is required to obtain a good approximate mapping that results in good regularization parameter choices.

Furthermore, the simplicity of our proposed method makes it widely applicable to many different fields of applications.
Future work includes extensions to learning parameters for hybrid projection methods or multi-parameter regularizers.  In addition, we plan to incorporate recent works on physics informed neural networks to design better networks and to develop methods to estimate the number of inclusions in addition to the regularity of inclusions in an image for further image analysis.

\bigskip
{\bf Acknowledgement:} This work was partially supported by the National Science Foundation (NSF) under grant DMS-1654175 (J. Chung), DMS-1723005 (M. Chung and J. Chung), and partially supported by The Villum Foundation grant no. 25893 (B. M. Afkham). This work was initiated as a part of the Statistical and Applied Mathematical Sciences Institute (SAMSI) Program on Numerical Analysis in Data Science in 2020 under the NSF grant DMS-1638521. Any opinions, findings, and conclusions or recommendations expressed in this material are those of the authors and do not necessarily reflect the views of the National Science Foundation.
\printbibliography[title=References]

@article{liu2021machine,
  title={Machine-learning-based prediction of regularization parameters for seismic inverse problems},
  author={Liu, Shihuan and Zhang, Jiashu},
  journal={Acta Geophysica},
  pages={1--12},
  year={2021},
  publisher={Springer}
}

@book{goodfellow2016machine,
  title={Machine Learning},
  author={Goodfellow, Ian and Bengio, Yoshua and Courville, Aaron},
%   journal={Deep learning},
%   volume={1},
%   pages={98--164},
  year={2016},
  publisher={MIT Press},
%   address = {Boston}
}

@article{robbins1951stochastic,
  title={A stochastic approximation method},
  author={Robbins, Herbert and Monro, Sutton},
  journal={The Annals of Mathematical Statistics},
  pages={400--407},
  year={1951},
  publisher={JSTOR},
  volume = {22},
  number = {3}
}

@article{kingma2014adam,
  title={Adam: A method for stochastic optimization},
  author={Kingma, Diederik P and Ba, Jimmy},
  journal={arXiv preprint arXiv:1412.6980},
  year={2014}
}

@article{kleywegt2002sample,
  title={The sample average approximation method for stochastic discrete optimization},
  author={Kleywegt, Anton J and Shapiro, Alexander and Homem-de-Mello, Tito},
  journal={SIAM Journal on Optimization},
  volume={12},
  number={2},
  pages={479--502},
  year={2002},
  publisher={SIAM}
}

@article{hornik1989multilayer,
  title={Multilayer feedforward networks are universal approximators},
  author={Hornik, Kurt and Stinchcombe, Maxwell and White, Halbert},
  journal={Neural Networks},
  volume={2},
  number={5},
  pages={359--366},
  year={1989},
  publisher={Elsevier}
}

@book{parker1994geophysical,
  title={Geophysical Inverse Theory},
  author={Parker, Robert L and Parker, Robert L},
  volume={1},
  year={1994},
  publisher={Princeton University Press}
}

@article{constable1987occam,
  title={Occam’s inversion: A practical algorithm for generating smooth models from electromagnetic sounding data},
  author={Constable, Steven C and Parker, Robert L and Constable, Catherine G},
  journal={Geophysics},
  volume={52},
  number={3},
  pages={289--300},
  year={1987},
  publisher={Society of Exploration Geophysicists}
}

@article{horesh2010optimal,
  title={Optimal experimental design for the large-scale nonlinear ill-posed problem of impedance imaging},
  author={Horesh, Lior and Haber, Eldad and Tenorio, Luis},
  journal={Large-Scale Inverse Problems and Quantification of Uncertainty},
  pages={273--290},
  year={2010},
  publisher={Wiley Online Library}
}

@article{gazzola2020inner,
  title={An inner--outer iterative method for edge preservation in image restoration and reconstruction},
  author={Gazzola, Silvia and Kilmer, Misha E and Nagy, James G and Semerci, Oguz and Miller, Eric L},
  journal={Inverse Problems},
  volume={36},
  number={12},
  pages={124004},
  year={2020},
  publisher={IOP Publishing}
}

@article{gazzola2014generalized,
  title={Generalized Arnoldi--Tikhonov method for sparse reconstruction},
  author={Gazzola, Silvia and Nagy, James G},
  journal={SIAM Journal on Scientific Computing},
  volume={36},
  number={2},
  pages={B225--B247},
  year={2014},
  publisher={SIAM}
}

@article{chung2019flexible,
  title={Flexible Krylov methods for $\ell_p$ regularization},
  author={Chung, Julianne and Gazzola, Silvia},
  journal={SIAM Journal on Scientific Computing},
  volume={41},
  number={5},
  pages={S149--S171},
  year={2019},
  publisher={SIAM}
}

@article{chung2008weighted,
  title={A weighted GCV method for Lanczos hybrid regularization},
  author={Chung, Julianne and Nagy, James G and O'Leary, Dianne P},
  journal={Electronic Transactions on Numerical Analysis},
  volume={28},
  number={149-167},
  pages={2008},
  year={2008}
}

@article{lin2010upre,
  title={UPRE method for total variation parameter selection},
  author={Lin, Youzuo and Wohlberg, Brendt and Guo, Hongbin},
  journal={Signal Processing},
  volume={90},
  number={8},
  pages={2546--2551},
  year={2010},
  publisher={Elsevier}
}

@article{wen2011parameter,
  title={Parameter selection for total-variation-based image restoration using discrepancy principle},
  author={Wen, You-Wei and Chan, Raymond H},
  journal={IEEE Transactions on Image Processing},
  volume={21},
  number={4},
  pages={1770--1781},
  year={2011},
  publisher={IEEE}
}

@article{langer2017automated,
  title={Automated parameter selection for total variation minimization in image restoration},
  author={Langer, Andreas},
  journal={Journal of Mathematical Imaging and Vision},
  volume={57},
  number={2},
  pages={239--268},
  year={2017},
  publisher={Springer}
}

@article{liao2009selection,
  title={Selection of regularization parameter in total variation image restoration},
  author={Liao, Haiyong and Li, Fang and Ng, Michael K},
  journal={JOSA A},
  volume={26},
  number={11},
  pages={2311--2320},
  year={2009},
  publisher={Optical Society of America}
}

@article{hammernik2018learning,
  title={Learning a variational network for reconstruction of accelerated MRI data},
  author={Hammernik, Kerstin and Klatzer, Teresa and Kobler, Erich and Recht, Michael P and Sodickson, Daniel K and Pock, Thomas and Knoll, Florian},
  journal={Magnetic Resonance in Medicine},
  volume={79},
  number={6},
  pages={3055--3071},
  year={2018},
  publisher={Wiley Online Library}
}

@inproceedings{zhang2017learning,
  title={Learning deep CNN denoiser prior for image restoration},
  author={Zhang, Kai and Zuo, Wangmeng and Gu, Shuhang and Zhang, Lei},
  booktitle={Proceedings of the IEEE Conference on Computer Vision and Pattern Recognition},
  pages={3929--3938},
  year={2017}
}

@article{lucas2018using,
  title={Using deep neural networks for inverse problems in imaging: beyond analytical methods},
  author={Lucas, Alice and Iliadis, Michael and Molina, Rafael and Katsaggelos, Aggelos K},
  journal={IEEE Signal Processing Magazine},
  volume={35},
  number={1},
  pages={20--36},
  year={2018},
  publisher={IEEE}
}

@article{mccann2017convolutional,
  title={Convolutional neural networks for inverse problems in imaging: A review},
  author={McCann, Michael T and Jin, Kyong Hwan and Unser, Michael},
  journal={IEEE Signal Processing Magazine},
  volume={34},
  number={6},
  pages={85--95},
  year={2017},
  publisher={IEEE}
}

@article{haber2012numerical,
  title={Numerical methods for A-optimal designs with a sparsity constraint for ill-posed inverse problems},
  author={Haber, Eldad and Magnant, Zhuojun and Lucero, Christian and Tenorio, Luis},
  journal={Computational Optimization and Applications},
  volume={52},
  number={1},
  pages={293--314},
  year={2012},
  publisher={Springer}
}

@article{calatroni2017bilevel,
  title={Bilevel approaches for learning of variational imaging models},
  author={Calatroni, Luca and Cao, Chung and De Los Reyes, Juan Carlos and Sch{\"o}nlieb, Carola-Bibiane and Valkonen, Tuomo},
  journal={Variational Methods: In Imaging and Geometric Control},
  volume={18},
  number={252},
  pages={2},
  year={2017},
  publisher={Walter de Gruyter GmbH}
}

@article{o1981bidiagonalization,
  title={A bidiagonalization-regularization procedure for large scale discretizations of ill-posed problems},
  author={O'Leary, Dianne P and Simmons, John A},
  journal={SIAM Journal on Scientific and Statistical Computing},
  volume={2},
  number={4},
  pages={474--489},
  year={1981},
  publisher={SIAM}
}

@article{bjorck1988bidiagonalization,
  title={A bidiagonalization algorithm for solving large and sparse ill-posed systems of linear equations},
  author={Bj{\"o}rck, {\AA}ke},
  journal={BIT Numerical Mathematics},
  volume={28},
  number={3},
  pages={659--670},
  year={1988},
  publisher={Springer}
}

@book{calvetti2007introduction,
  title={An Introduction to Bayesian Scientific Computing: Ten Lectures on Subjective Computing},
  author={Calvetti, Daniela and Somersalo, Erkki},
  volume={2},
  year={2007},
  publisher={Springer Science \& Business Media},
%   address = {New York}
}

@article{antil2020bilevel,
  title={Bilevel optimization, deep learning and fractional Laplacian regularization with applications in tomography},
  author={Antil, Harbir and Di, Zichao Wendy and Khatri, Ratna},
  journal={Inverse Problems},
  volume={36},
  number={6},
  pages={064001},
  year={2020},
  publisher={IOP Publishing}
}

@book{shapiro2014lectures,
  title={Lectures on Stochastic Programming: Modeling and Theory},
  author={Shapiro, Alexander and Dentcheva, Darinka and Ruszczy{\'n}ski, Andrzej},
  year={2014},
  publisher={SIAM},
%   address = {Philadelphia}
}

@article{cybenko1989approximation,
  title={Approximation by superpositions of a sigmoidal function},
  author={Cybenko, George},
  journal={Mathematics of Control, Signals and Systems},
  volume={2},
  number={4},
  pages={303--314},
  year={1989},
  publisher={Springer}
}

@article{farquharson2004comparison,
  title={A comparison of automatic techniques for estimating the regularization parameter in non-linear inverse problems},
  author={Farquharson, Colin G and Oldenburg, Douglas W},
  journal={Geophysical Journal International},
  volume={156},
  number={3},
  pages={411--425},
  year={2004},
  publisher={Blackwell Science Ltd Oxford, UK}
}

@article{vogel1996non,
  title={Non-convergence of the L-curve regularization parameter selection method},
  author={Vogel, Curtis R},
  journal={Inverse Problems},
  volume={12},
  number={4},
  pages={535},
  year={1996},
  publisher={IOP Publishing}
}

@book{bardsley2018computational,
  title={Computational Uncertainty Quantification for Inverse Problems},
  author={Bardsley, Johnathan M},
%   volume={19},
  year={2018},
%   address = {Philadelphia},
  publisher={SIAM}
}

@article{donoho1995noising,
  title={De-noising by soft-thresholding},
  author={Donoho, David L},
  journal={IEEE Transactions on Information Theory},
  volume={41},
  number={3},
  pages={613--627},
  year={1995},
  publisher={IEEE}
}

@article{luiken2020comparing,
  title={Comparing RSVD and Krylov methods for linear inverse problems},
  author={Luiken, Nick and van Leeuwen, Tristan},
  journal={Computers \& Geosciences},
  volume={137},
  pages={104427},
  year={2020},
  publisher={Elsevier}
}

@article{arridge2019solving,
  title={Solving inverse problems using data-driven models},
  author={Arridge, Simon and Maass, Peter and {\"O}ktem, Ozan and Sch{\"o}nlieb, Carola-Bibiane},
  journal={Acta Numerica},
  volume={28},
  pages={1--174},
  year={2019},
  publisher={Cambridge University Press}
}

@article{chung2011designing,
  title={Designing optimal spectral filters for inverse problems},
  author={Chung, Julianne and Chung, Matthias and O'Leary, Dianne P},
  journal={SIAM Journal on Scientific Computing},
  volume={33},
  number={6},
  pages={3132--3152},
  year={2011},
  publisher={SIAM}
}

@article{de2017bilevel,
  title={Bilevel parameter learning for higher-order total variation regularisation models},
  author={De los Reyes, Juan Carlos and Sch{\"o}nlieb, C-B and Valkonen, Tuomo},
  journal={Journal of Mathematical Imaging and Vision},
  volume={57},
  number={1},
  pages={1--25},
  year={2017},
  publisher={Springer}
}

@Book{Pukelsheim2006,
  Title                    = {Optimal Design of Experiments},
  Author                   = {F Pukelsheim},
  Publisher                = {SIAM},
  Year                     = {2006},
%   Address                  = {Philadelphia},
  Owner                    = {mcchung},
  Pages                    = {454},
  Timestamp                = {2012.12.01}
}

@Book{Atkinson2007,
  Title                    = {Optimum Experimental Designs, with SAS},
  Author                   = {A~C Atkinson and A~N Donev and R Tobias},
  Publisher                = {Oxford University Press},
  Year                     = {2007},
%   Address                  = {Oxford},
  Owner                    = {mcchung},
  Pages                    = {528},
  Timestamp                = {2012.12.01}
}

@article{wang2020learning,
  title={Learning priors for adversarial autoencoders},
  author={Wang, Hui-Po and Peng, Wen-Hsiao and Ko, Wei-Jan},
  journal={APSIPA Transactions on Signal and Information Processing},
  volume={9},
  year={2020},
  publisher={Cambridge University Press}
}

@book{tenorio2017introduction,
  title={An Introduction to Data Analysis and Uncertainty Quantification for Inverse Problems},
  author={Tenorio, Luis},
  year={2017},
  publisher={SIAM}
}

@article{mead2008newton,
  title={A Newton root-finding algorithm for estimating the regularization parameter for solving ill-conditioned least squares problems},
  author={Mead, Jodi L and Renaut, Rosemary A},
  journal={Inverse Problems},
  volume={25},
  number={2},
  pages={025002},
  year={2008},
  publisher={IOP Publishing}
}

@article{galatsanos1992methods,
  title={Methods for choosing the regularization parameter and estimating the noise variance in image restoration and their relation},
  author={Galatsanos, Nikolas P and Katsaggelos, Aggelos K},
  journal={IEEE Transactions on Image Processing},
  volume={1},
  number={3},
  pages={322--336},
  year={1992}
}

@book{engl1996regularization,
  title={Regularization of Inverse Problems},
  author={Engl, Heinz Werner and Hanke, Martin and Neubauer, Andreas},
%   volume={375},
  year={1996},
  publisher={Springer Science \& Business Media},
%   address = {New York}
}

@book{hansen2010discrete,
  title={Discrete Inverse Problems: Insight and Algorithms},
  author={Hansen, Per Christian},
  year={2010},
  publisher={SIAM}
}

@book{hansen2006deblurring,
  title={Deblurring Images: Matrices, Spectra, and Filtering},
  author={Hansen, P.C. and Nagy, J.G. and O'Leary, D.P.},
%   isbn={9780898716184},
%   lccn={2006044387},
%   series={Fundamentals of Algorithms},
%   url={https://books.google.dk/books?id=JjbSoRR9T-0C},
  year={2006},
  publisher={SIAM},
%   address = {Philadelphia}
}

@article{rudin1992nonlinear,
  title={Nonlinear total variation based noise removal algorithms},
  author={Rudin, Leonid I and Osher, Stanley and Fatemi, Emad},
  journal={Physica D: Nonlinear Phenomena},
  volume={60},
  number={1-4},
  pages={259--268},
  year={1992},
  publisher={North-Holland}
}

@incollection{lamm2000survey,
  title={A survey of regularization methods for first-kind Volterra equations},
  author={Lamm, Patricia K},
  booktitle={Surveys on Solution Methods for Inverse Problems},
  pages={53--82},
  year={2000},
  publisher={Springer}
}

@online{regTools,
  author = {Per Christian Hansen},
  title = {Regtools},
  year = 2020,
  url = {https://www.mathworks.com/matlabcentral/fileexchange/52-regtools},
  urldate = {2020-11-08}
}

@online{randomSheppLogan,
  author = {Matthias Chung},
  title = {Random-Shepp-Logan-Phantom},
  year = 2020,
  url = {https://github.com/matthiaschung/Random-Shepp-Logan-Phantom},
  urldate = {2020-12-14}
}

@article{ruthotto2018optimal,
  title={Optimal experimental design for inverse problems with state constraints},
  author={Ruthotto, Lars and Chung, Julianne and Chung, Matthias},
  journal={SIAM Journal on Scientific Computing},
  volume={40},
  number={4},
  pages={B1080--B1100},
  year={2018},
  publisher={SIAM}
}

@book{gramacy2020surrogates,
  title={Surrogates: Gaussian Process Modeling, Design, and Optimization for the Applied Sciences},
  author={Gramacy, Robert B},
  year={2020},
  publisher={CRC Press}
}

@article{haber2000gcv,
  title={A GCV based method for nonlinear ill-posed problems},
  author={Haber, Eldad and Oldenburg, Douglas},
  journal={Computational Geosciences},
  volume={4},
  number={1},
  pages={41--63},
  year={2000},
  publisher={Springer}
}

@article{de2016machine,
  title={A machine learning approach to optimal Tikhonov regularisation I: Affine manifolds},
  author={De Vito, Ernesto and Fornasier, Massimo and Naumova, Valeriya},
  journal={arXiv preprint arXiv:1610.01952},
  year={2016}
}

@article{li2020nett,
  title={NETT: Solving inverse problems with deep neural networks},
  author={Li, Housen and Schwab, Johannes and Antholzer, Stephan and Haltmeier, Markus},
  journal={Inverse Problems},
  volume={36},
  number={6},
  pages={065005},
  year={2020},
  publisher={IOP Publishing}
}

@article{haber2003learning,
  title={Learning regularization functionals-a supervised training approach},
  author={Haber, E and Tenorio, L},
  journal={Inverse Problems},
  volume={19},
  number={3},
  pages={611},
  year={2003},
  publisher={IOP Publishing}
}

@article{haber2009numerical,
  title={Numerical methods for the design of large-scale nonlinear discrete ill-posed inverse problems},
  author={Haber, Eldad and Horesh, Lior and Tenorio, Luis},
  journal={Inverse Problems},
  volume={26},
  number={2},
  pages={025002},
  year={2009},
  publisher={IOP Publishing}
}

@article{haber2008numerical,
  title={Numerical methods for experimental design of large-scale linear ill-posed inverse problems},
  author={Haber, Eldad and Horesh, Lior and Tenorio, Luis},
  journal={Inverse Problems},
  volume={24},
  number={5},
  pages={055012},
  year={2008},
  publisher={IOP Publishing}
}

@article{chung2012optimal,
  title={Optimal filters from calibration data for image deconvolution with data acquisition error},
  author={Chung, Julianne and Chung, Matthias and O'Leary, Dianne P},
  journal={Journal of Mathematical Imaging and Vision},
  volume={44},
  number={3},
  pages={366--374},
  year={2012},
  publisher={Springer}
}

@ARTICLE{5701777,
  author={W. {Dong} and L. {Zhang} and G. {Shi} and X. {Wu}},
  journal={IEEE Transactions on Image Processing}, 
  title={Image deblurring and super-resolution by adaptive sparse domain selection and adaptive regularization}, 
  year={2011},
  volume={20},
  number={7},
  pages={1838-1857},
%   doi={10.1109/TIP.2011.2108306}
}

@article{Huang2008,
%   doi = {10.1137/070703533},
%   url = {https://doi.org/10.1137/070703533},
  year = {2008},
%   month = jan,
  publisher = {Society for Industrial {\&} Applied Mathematics ({SIAM})},
  volume = {7},
  number = {2},
  pages = {774--795},
  author = {Yumei Huang and Michael K. Ng and You-Wei Wen},
  title = {A Fast Total Variation Minimization Method for Image Restoration},
  journal = {Multiscale Modeling {\&} Simulation}
}

@ARTICLE{7448477,
  author={J. {Pan} and Z. {Hu} and Z. {Su} and M. {Yang}},
  journal={IEEE Transactions on Pattern Analysis and Machine Intelligence}, 
  title={$L_0$-regularized intensity and gradient prior for deblurring text images and beyond}, 
  year={2017},
  volume={39},
  number={2},
  pages={342-355},
%   doi={10.1109/TPAMI.2016.2551244}
}

@INPROCEEDINGS{6909767,
  author={J. {Pan} and Z. {Hu} and Z. {Su} and M. {Yang}},
  booktitle={2014 IEEE Conference on Computer Vision and Pattern Recognition}, 
  title={Deblurring text images via $L_0$-regularized intensity and gradient prior}, 
  year={2014},
  volume={},
  number={},
  pages={2901-2908},
%   doi={10.1109/CVPR.2014.371}
}

@article{Pearson2018,
%   doi = {10.1051/emsci/2017010},
%   url = {https://doi.org/10.1051/emsci/2017010},
  year = {2018},
  publisher = {{EDP} Sciences},
  volume = {2},
  pages = {1},
  author = {James Pearson and Clara Pennock and Tom Robinson},
  title = {Auto-detection of strong gravitational lenses using convolutional neural networks},
  journal = {Emergent Scientist}
}

@article{1930-8337_2016_4_1007,
title = {The Bayesian formulation of EIT: Analysis and algorithms},
journal = {Inverse Problems \& Imaging},
volume = {10},
number = {4},
pages = {1007},
year = {2016},
% issn = {1930-8337},
% doi = {10.3934/ipi.2016030},
% url = {http://aimsciences.org//article/id/d9e2dc77-0da0-4532-95df-b2c6e5a9eece},
author = {Matthew M.  Dunlop and Andrew M.  Stuart}
}

@article{BuiThanh2014,
%   doi = {10.1137/120894877},
%   url = {https://doi.org/10.1137/120894877},
  year = {2014},
%   month = jan,
  publisher = {Society for Industrial {\&} Applied Mathematics ({SIAM})},
  volume = {2},
  number = {1},
  pages = {203--222},
  author = {Tan Bui-Thanh and Omar Ghattas},
  title = {An analysis of infinite dimensional Bayesian inverse shape acoustic scattering and its numerical approximation},
  journal = {{SIAM}/{ASA} Journal on Uncertainty Quantification}
}

@book{ibragimov2012gaussian,
  title={Gaussian Random Processes},
  author={Ibragimov, Ildar Abdulovich and Rozanov, Yurii Antol'evich},
  volume={9},
  year={2012},
  publisher={Springer Science \& Business Media},
%   address ={New York}
}

@article{huang2006extreme,
  title={Extreme learning machine: theory and applications},
  author={Huang, Guang-Bin and Zhu, Qin-Yu and Siew, Chee-Kheong},
  journal={Neurocomputing},
  volume={70},
  number={1-3},
  pages={489--501},
  year={2006},
  publisher={Elsevier}
}

@article{paige1982lsqr,
  title={LSQR: An algorithm for sparse linear equations and sparse least squares},
  author={Paige, Christopher C and Saunders, Michael A},
  journal={ACM Transactions on Mathematical Software (TOMS)},
  volume={8},
  number={1},
  pages={43--71},
  year={1982},
  publisher={ACM New York, NY, USA}
}

@article{gazzola2019ir,
  title={{IR Tools}: a {MATLAB} package of iterative regularization methods and large-scale test problems},
  author={Gazzola, Silvia and Hansen, Per Christian and Nagy, James G},
  journal={Numerical Algorithms},
  volume={81},
  number={3},
  pages={773--811},
  year={2019},
  publisher={Springer}
}

@article{min2013inverse,
  title={Inverse estimation of the initial condition for the heat equation},
  author={Min, Tao and Geng, Bei and Ren, Jucheng},
  journal={International Journal of Pure and Applied Mathematics},
  volume={82},
  number={4},
  pages={581--593},
  year={2013},
  publisher={Academic Publications, Ltd.}
}

@book{natterer2001mathematics,
  title={The Mathematics of Computerized Tomography},
  author={Natterer, Frank},
  year={2001},
  publisher={SIAM},
%   address = {Philadelphia}
}

@article{borsic2009vivo,
  title={In vivo impedance imaging with total variation regularization},
  journal={IEEE Transactions on Medical Imaging},
  volume={29},
  number={1},
  pages={44--54},
  year={2009},
  publisher={IEEE}
}

@book{kak2002principles,
  title={Principles of Computerized Tomographic Imaging},
  author={Kak, Avinash C and Slaney, Malcolm and Wang, Ge},
  year={2002},
  publisher={Wiley Online Library},
%   address = {Hoboken}
}

@article{hansen2018air,
  title={{AIR} {T}ools {II}: algebraic iterative reconstruction methods, improved implementation},
  author={Hansen, Per Christian and J{\o}rgensen, Jakob Sauer},
  journal={Numerical Algorithms},
  volume={79},
  number={1},
  pages={107--137},
  year={2018},
  publisher={Springer}
}

@article{dittmer2020regularization,
  title={Regularization by architecture: A deep prior approach for inverse problems},
  author={Dittmer, S{\"o}ren and Kluth, Tobias and Maass, Peter and Baguer, Daniel Otero},
  journal={Journal of Mathematical Imaging and Vision},
  volume={62},
  number={3},
  pages={456--470},
  year={2020},
  publisher={Springer}
}

@article{hornik1990universal,
  title={Universal approximation of an unknown mapping and its derivatives using multilayer feedforward networks},
  author={Hornik, Kurt and Stinchcombe, Maxwell and White, Halbert},
  journal={Neural Networks},
  volume={3},
  number={5},
  pages={551--560},
  year={1990},
  publisher={Elsevier}
}

@article{Puetter2005,
%   doi = {10.1146/annurev.astro.43.112904.104850},
%   url = {https://doi.org/10.1146/annurev.astro.43.112904.104850},
  year = {2005},
%   month = sep,
  publisher = {Annual Reviews},
  volume = {43},
  number = {1},
  pages = {139--194},
  author = {R.C. Puetter and T.R. Gosnell and Amos Yahil},
  title = {Digital image reconstruction: Deblurring and denoising},
  journal = {Annual Review of Astronomy and Astrophysics}
}

\end{document}